\documentclass[manuscript,screen]{acmart}

\usepackage{graphicx}
\usepackage{subcaption}
\usepackage{listings}
\usepackage{soul}
\graphicspath{{figdata/}}

\definecolor{codegreen}{rgb}{0,0.6,0}
\definecolor{codegray}{rgb}{0.5,0.5,0.5}
\definecolor{codepurple}{rgb}{0.58,0,0.82}
\definecolor{backcolour}{rgb}{0.95,0.95,0.92}

\lstdefinestyle{mystyle}{
    backgroundcolor=\color{backcolour},   
    commentstyle=\color{codegreen},
    keywordstyle=\color{magenta},
    numberstyle=\tiny\color{codegray},
    stringstyle=\color{codepurple},
    basicstyle=\ttfamily\footnotesize,
    breakatwhitespace=false,         
    breaklines=true,                 
    captionpos=b,                    
    keepspaces=true,                 
    numbers=left,                    
    numbersep=5pt,                  
    showspaces=false,                
    showstringspaces=false,
    showtabs=false,                  
    tabsize=2
}
\AtBeginDocument{%
  }

\setcopyright{acmcopyright}
\copyrightyear{2023}
\acmYear{2023}
\acmDOI{XXXXXXX.XXXXXXX}

\begin{document}

\title{Algorithm xxxx: HiPPIS A High-Order Positivity-Preserving Mapping Software for Structured Meshes}

\author{Timbwoga A. J. Ouermi}
\email{touermi@cs.utah.edu}
\orcid{1234-5678-9012}
\author{Robert M. Kirby}
\email{kirby@cs.utah.edu}
\author{Martin Berzins}
\email{mb@sci.utah.edu}
\affiliation{%
  \institution{University of Utah Scientific Computing Imaging Institute}
  \streetaddress{72 Central campus Drive}
  \city{Salt Lake City}
  \state{Utah}
  \country{USA}
  \postcode{84112}
}

%\author{Lars Th{\o}rv{\"a}ld}
%\affiliation{%
%  \institution{The Th{\o}rv{\"a}ld Group}
%  \streetaddress{1 Th{\o}rv{\"a}ld Circle}
%\email{larst@affiliation.org}
%\affiliation{%
%  \institution{Inria Paris-Rocquencourt}
%\affiliation{%
% \institution{Rajiv Gandhi University}
%\affiliation{%
%  \institution{Tsinghua University}
%  \streetaddress{30 Shuangqing Rd}
%\affiliation{%
%  \institution{Palmer Research Laboratories}
%  \streetaddress{8600 Datapoint Drive}
%  \country{USA}
%\affiliation{%
%  \institution{The Th{\o}rv{\"a}ld Group}
%  \streetaddress{1 Th{\o}rv{\"a}ld Circle}
%\email{jsmith@affiliation.org}
%\affiliation{%
%  \institution{The Kumquat Consortium}
%  \city{New York}
%\email{jpkumquat@consortium.net}

\renewcommand{\shortauthors}{Ouermi et al.}

\begin{abstract}
Polynomial interpolation is an important component of many computational problems. 
In several of these computational problems, failure to preserve positivity when using polynomials to approximate or map data values between meshes can lead to negative unphysical quantities. 
Currently, most polynomial-based methods for enforcing positivity are based on splines and polynomial rescaling. 
The spline-based approaches build interpolants that are positive over the intervals in which they are defined and may require solving a minimization problem and/or system of equations. 
The linear polynomial rescaling methods allow for high-degree polynomials but enforce positivity only at limited locations (e.g., quadrature nodes). 
This work introduces open-source software (HiPPIS) for high-order data-bounded interpolation (DBI) and positivity-preserving interpolation (PPI) that addresses the limitations of both the spline and polynomial rescaling methods. 
HiPPIS is suitable for approximating and mapping physical quantities such as mass, density, and concentration between meshes while preserving positivity. 
This work provides Fortran and Matlab implementations of the DBI and PPI methods, presents an analysis of the mapping error in the context of PDEs, and uses several 1D and 2D numerical examples to demonstrate the benefits and limitations of HiPPIS. %\color{red}\st{, and introduces different strategies to improve locality, vectorization, and overall performance.}\color{black} 
\end{abstract}

\begin{CCSXML}
<ccs2012>
   <concept>
       <concept_id>10002950.10003714.10003715.10003720</concept_id>
       <concept_desc>Mathematics of computing~Computations on polynomials</concept_desc>
       <concept_significance>500</concept_significance>
       </concept>
 </ccs2012>
\end{CCSXML}

\ccsdesc[500]{Mathematics of computing~Computations on polynomials}

%\begin{CCSXML}
%<ccs2012>
%  <concept_id>10010520.10010553.10010562</concept_id>
%  <concept_desc>Computer systems organization~Embedded systems</concept_desc>
%  <concept_significance>500</concept_significance>
%  <concept_id>10010520.10010575.10010755</concept_id>
%  <concept_desc>Computer systems organization~Redundancy</concept_desc>
%  <concept_significance>300</concept_significance>
%  <concept_id>10010520.10010553.10010554</concept_id>
%  <concept_desc>Computer systems organization~Robotics</concept_desc>
%  <concept_significance>100</concept_significance>
%  <concept_id>10003033.10003083.10003095</concept_id>
%  <concept_desc>Networks~Network reliability</concept_desc>
%  <concept_significance>100</concept_significance>
%</ccs2012>
%\end{CCSXML}
%
%\ccsdesc[500]{Computer systems organization~Embedded systems}

\keywords{positivity-preserving, data-bounded, polynomial interpolation, vectorization}

\maketitle

\section{Introduction}
\label{sec:intro}
  Mapping data values from one grid to another is a fundamental part of many computational problems. 
  Preserving certain properties such as positivity when interpolating solution values between meshes is important. 
  In many applications ~\cite{skamrock, LIU201982, SUBBAREDDY2017827, BALSARA20127504, Damodar, ZHANG2017}, failure to preserve the positivity of quantities such as mass, density, and concentration results in negative values that are unphysical.
  These negative values may propagate to other calculations and corrupt other quantities. 
  Many polynomial-based methods have been developed to address these limitations.  
  
  %% Linear rescaling methods
  Positivity-preserving methods based on linear polynomial rescaling are introduced in ~\cite{Light, ZHANG2017, Zhang2012_2, LIU201982, HU2013169}.
  These polynomial rescaling methods are often used in the context of hyperbolic PDEs, in numerical weather prediction (NWP) ~\cite{Light}, combustion simulation ~\cite{LIU201982, HU2013169}, and other applications.   
  These methods introduce rescaling parameters obtained from quadrature weights that are used to linearly rescale the polynomial to ensure positivity at the quadrature nodes and conserve mass. 
  These approaches ensure positivity only at the set of mesh points used for the simulation but do not address the case of mapping data values between different meshes, which is the focus of HiPPIS. 

  %% Spline based approaches
  Other approaches for preserving positivity that are based on splines can be found in computer-aided design (CAD), graphics, and visualization ~\cite{Schmidt1988, Schmidt1987, HUSSAIN2008446, karim2015positivity, SARFRAZ199269, fritsch1980monotone}. 
  Several positivity- and monotonicity-preserving cubic splines have been developed. 
  %
  %These cubic spline approaches provide a sufficient and/or necessary condition for monotonicity, positivity, and $\mathcal{C}^{1}$ continuity. 
  %
  A widely used example of such an approach is the piecewise cubic Hermite interpolation (PCHIP) ~\cite{fritsch1980monotone}, which is available as open-source code in ~\cite{moler2004numerical}. 
  In addition, quartic and quintic spline-based approaches have been introduced in ~\cite{HE199451, Lux2019ANAF, Lux2023, hussain2009c2, Hussain2018ACR}. 
  These methods impose some restrictions on the first and second derivatives to ensure monotonicity, positivity, and continuity.
  For instance, the monotonic quintic spline interpolation (MQSI) methods in ~\cite{Lux2019ANAF} and \cite{Lux2023} use the sufficient conditions stated in ~\cite{Schmidt1988} and ~\cite{doi:10.1137/0915035} to check for monotonicity and iteratively adjust the first and second derivative values to enforce monotonicity.

  Positivity can also be enforced using ENO-type methods ~\cite{Damodar, berzins2010nonlinear, Berzins, Ouermi2022}, which enforce data-boundedness and positivity by adaptively selecting mesh points to build the stencil used to construct the positive interpolant for each interval.  
  ENO-type methods use divided differences to develop a sufficient condition for data-boundedness or positivity that is used to guide the stencil selection process. 
  The software introduced in this work is based on the high-order ENO-type data-bounded interpolation (DBI) and positivity-preserving interpolation (PPI) methods in ~\cite{Ouermi2022}.
  The work in ~\cite{Ouermi2022} provides a positivity-preserving method that uses higher degree polynomials compared to the other ENO-type methods in ~\cite{Damodar, berzins2010nonlinear, Berzins} and the spline-based methods.

  Given that polynomial interpolation is well-established and widely used, several implementations of the different polynomial approximation algorithms are available.
  For example, FunC by Green et al. \cite{Green2019} uses polynomial interpolation with a lookup table for faster approximation of a given function compared to direct evaluation.
  However, most of these implementations, including FunC, do not preserve data-boundedness and positivity.  
  The implementations available for positivity preservation are based on splines ~\cite{HE199451, fritsch1980monotone} and polynomial rescaling ~\cite{ZHANG2017, Light}.
  The spline-based approaches often require solving a linear system of equations to ensure continuity and an optimization problem in the case of quartic and quintic splines. 
  These spline approaches are often limited to fifth-order polynomials and can be computationally expensive in cases where solving a global optimization problem is required. 
  A full suite of test problems comparing the DBI and PPI methods against different spline-based methods including PCHIP ~\cite{fritsch1980monotone}, MQSI ~\cite{Lux2019ANAF}, and shape-preserving splines (SPS) ~\cite{10.1145/264029.264059} has been undertaken by the authors in ~\cite{arkivtajo}. 
  The different polynomial rescaling methods allow for polynomial degrees higher than five and are built as part of larger partial differential equation (PDE) solvers ~\cite{Light, ZHANG2017}. 
  As previously mentioned, the polynomial rescaling approaches guarantee positivity only at a given set of points, not over the entire domain. 
  The present work provides an implementation of a high-order software (HiPPIS) based on ~\cite{Ouermi2022} that guarantees positivity over the entire domain where the interpolant is defined. 
  In addition, this work evaluates the use of HiPPIS in the context of function approximation and mapping between different meshes. 
  This evaluation provides an analysis of the mapping error in the case of PDEs and numerical examples demonstrating the benefits and limitations of HiPPIS.  

  The remaining parts of the paper are organized as follows: 
  Section \ref{sec:background} presents the background for the mathematical framework required for the DBI and PPI methods.
  Section \ref{sec:algorithm} provides the algorithms used to build the software, and the descriptions of the different components of HiPPIS. % \color{red} \st{, and the techniques used to enable vectorization, increase locality, and improve overall computational performance.} \color{black}
  Section \ref{sec:numerical-examples} provides several 1D and 2D numerical examples, while Section \ref{sec:error-mapping-examples} conducts an analysis and evaluation of the mapping error in the context of time-dependent PDEs. 
  A discussion and concluding remarks are presented in Section \ref{sec:discussion-conlusion}.

  \section{Mathematical Framework}
\label{sec:background}
  This section provides a summary and the theoretical background of both the DBI and PPI methods.
  \subsection{Adaptive Polynomial Construction}
    Both the DBI and PPI methods rely on the Newton polynomial \cite{doi:10.1137/0912034,10.2307/2004888} representation to build interpolants that are positive or bounded by the data values. 
    The ability to adaptively choose stencil points to construct the interpolation, as in ENO methods \cite{HARTEN19973}, is the key feature employed to develop the data-bounded and positivity-preserving interpolants. 

    Consider a 1D mesh defined as follows: 
    \begin{equation}\label{eq:RandomMesh}
       \mathcal{M} = \{ x_{i-J}, \cdots, x_{i}, x_{i+1}, \cdots, x_{i+L}\},
    \end{equation} 
    where $x_{i-J} < \cdots < x_{i} < x_{i+1}< \cdots < x_{i+L}$, and $\{ u_{i-J}, \cdots, u_{i+L} \}$ is the set of data values associated with the mesh points in Equation (\ref{eq:RandomMesh}).
    The subscripts $J$, $L$, $i, \in \mathbb{N}_{0} = \mathbb{N} \cup \{0\}$, and $x_{k}$, $u_{k} \in \mathbb{R}$ for $i-J \leq k \leq i+L$. 
    The DBI and PPI procedure starts by setting the initial stencil $\mathcal{V}_{0}$,
    \begin{equation}\label{eq:v0}
      \mathcal{V}_{0} = \{ x_{i}, x_{i+1}\} =\{ x_{0}^{l}, x_{0}^{r}\}. 
    \end{equation}
    The stencil $\mathcal{V}_{0}$ in Equation (\ref{eq:v0}) is expanded by successively appending a point to the right or left of $\mathcal{V}_{j}$ to form $\mathcal{V}_{j+1}$. 
    Once the final stencil $\mathcal{V}_{n-1}$ is obtained, the interpolant of degree $n$ defined on $I_{i}=\{x_{i}, x_{i+1}\}$ can be written as
    \begin{equation}\label{eq:newtonPoly}
      \begin{gathered}
      U_{n}(x)= \quad u_{i} + U[x_{0}^{l}, x_{0}^{r}] \pi_{0,i}(x) + 
            U[x_{1}^{l}, \cdots, x_{1}^{r}] \pi_{1,i}(x) +  
            \cdots +  U[x_{n -1}^{l}, \cdots, x_{n -1}^{r}]\pi_{n-1,i}(x),
      \end{gathered}
    \end{equation}
    where 
    $ \pi_{0,i}(x) = (x-x_{i}), \pi_{1,i}(x)=(x-x_{i})(x-x_{1}^{e}), \cdots$ are the Newton basis functions.
    $x_{j}^{e}$ is the point added to expand the stencil $\mathcal{V}_{j-2}$ to $\mathcal{V}_{j-1}$ and can be explicitly expressed as  
    \begin{equation*}\label{eq:xe}
       \begin{cases}
         x_{0}^{e} = x_{i}, \\%\quad
         x_{1}^{e} = x_{i+1}, \\% \quad
         x_{j}^{e} = \mathcal{V}_{j-1} \setminus \mathcal{V}_{j-2}, \quad 2 \leq j \leq n-1.
       \end{cases}
    \end{equation*}
    The divided differences are recursively defined as follows:
    \begin{equation*}\label{eq:divDiff}
      \begin{cases}
      U[x_{i}] = u_{i}\\
      U[x_{i}, \cdots, x_{i+j}] = \frac{U[x_{i+1}, \cdots,x_{i+j}]-U[x_{i}, \cdots,x_{i+j-1}]}{x_{i+j} -x_{i}}.
      \end{cases}
    \end{equation*} 
    The polynomial $U_{n}(x)$ can be compactly expressed as  
    \begin{equation}\label{eq:UN3}
      \begin{gathered}
      U_{n}(x) = u_{i} + (u_{i+1}-u_{i})S_{n}(x).
      \end{gathered}
    \end{equation} 
    $S_{n}(x)$ in Equation (\ref{eq:UN3}) is defined as
    \begin{equation}\label{eq:SN1}
      \begin{gathered}
      S_{n}(x) = s
                \bigg(1 + \frac{(s-1)}{d_{1}}\lambda_{1} 
                \bigg(1 + \frac{(s-t_{2})}{d_{2}}\lambda_{2} \bigg( \cdots %\lambda_{n-2} 
                \bigg(1 + \frac{(s-t_{n-1})}{d_{n-1}}\lambda_{n-1} \bigg) \cdots \bigg), 
      \end{gathered}
    \end{equation} 
    where $s$, $t_{j}$, and $d_{j}$ are expressed as follows:
    \begin{equation}\label{eq:s}
      0 \leq s = \frac{x-x_{i}}{x_{i+1}-x_{i}} = \frac{x-x_{0}^{e}}{x_{0}^{r}-x_{0}^{l}} \leq 1,
    \end{equation}
    \begin{equation}\label{eq:tk}
      t_{j} = -\frac{x_{i}-x_{j}^{e}}{x_{0}^{r}-x_{0}^{l}}, \textrm{ and } 
    \end{equation}
    \begin{equation}\label{eq:dk}
      0 \leq d_{j} = \frac{x_{j}^{r}-x_{j}^{l}}{x_{0}^{r}-x_{0}^{l}} .
    \end{equation}
    $s$ and $d_{j}$ in Equations (\ref{eq:SN1}), (\ref{eq:s}), and (\ref{eq:dk}) are defined such that $s \in [0, 1]$ and $d_{j} \geq 0$.
    %
    %In the results, we consider the term $\bar{\lambda}_{j}$, defined as
    The positivity-preserving and data-bounded interpolants are obtained by imposing some bounds on $\bar{\lambda}_{j}$, defined as
    \begin{equation} \label{eq:bar_lambda}
      \bar{\lambda}_{j} = \prod_{k=1}^{j} \lambda_{j} = \lambda_j \bar{\lambda}_{j-1} = \prod_{k=1}^{j}\lambda_{k} =
      \begin{cases}
        1 \quad j = 0\\
      \frac{U[x_{j}^{l}, \cdots, x_{j}^{r}]}{U[x_{0}^{l}, x_{0}^{r}]}\prod_{k=1}^{j}(x_{k}^{r}-x_{k}^{l}), \quad 1\leq j \leq n-1.
      \end{cases}
    \end{equation}

  \subsection{Positivity-Preserving and Data-Bounded Interpolation}
  \label{subsec:ppi-dbi}
  For a given interval inside a mesh, the DBI and PPI polynomial interpolant is constructed by adaptively selecting points near the target interval to build an interpolation stencil and polynomial that together meet the requirements for positivity or data-boundedness. 
  Requiring positivity alone can lead to large oscillations and extrema that degrade the approximation. 
  Positivity alone does not restrict how much the interpolant is allowed to grow beyond the data values. The large oscillations can be removed with the PCHIP, MQSI, and DBI methods. 
  However, in the case where a given interval $I_{i}$ has a hidden extremum, PCHIP, MQSI, and DBI will truncate the extremum.
  As in ~\cite{berzins2010nonlinear,sekora2009extremumpreserving}, the interval $I_{i}$ has an extremum when two  of the three divided differences  $\sigma_{i-1}=U[x_{i-1}, x_{i}]$, $\sigma_{i} =U[x_{i}, x_{i+1}]$, and $\sigma_{i+1}=U[x_{i+1}, x_{i+2}]$ of neighboring intervals are of opposite signs. 
  The constrained PPI algorithm addresses these limitations by allowing the constructed interpolant to grow beyond the data values but not produce extrema that are too large. 

  The positive polynomial interpolant is constrained as follows: 
  \begin{equation}\label{eq:boundUP}
    u_{min} \leq U^{p}(x)= u_{i} + (u_{i+1}-u_{i})S_{n}(x) \leq u_{max}.
  \end{equation} 
  The bounds $u_{min}$ and $u_{max}$ in Equation (\ref{eq:boundUP}) are defined as 
  \begin{equation}\label{eq:uminumax}
    \begin{cases}
       u_{min} =  min ( u_{i}, u_{i+1} ) -\Delta_{min}, \\ 
       u_{max} =  max (u_{i}, u_{i+1} ) + \Delta_{max}.
     \end{cases}
  \end{equation}
  The  parameters $\Delta_{min}$ and $\Delta_{max}$ in Equation (\ref{eq:uminumax}) are positive, and the data-bounded interpolant is obtained for 
 $\Delta_{min}$=$\Delta_{max} =0.0$.
 These parameters are chosen according to
  \begin{equation} \label{eq:umin2}
    \Delta_{min} =
    \begin{cases}
      \epsilon_{1}\big|min\big(u_{i}, u_{i+1}\big) \big| & \textrm{if } \sigma_{i-1} \sigma_{i+1} < 0 
                                               \textrm{ and } \sigma_{i-1} < 0 
                                               \textrm{ or } \sigma_{i-1}\sigma_{i+1} \geq 0 
                                               \textrm{ and } \sigma_{i-1} \sigma_{i} < 0 \\
      \epsilon_{0} \big|min\big(u_{i}, u_{i+1}\big) \big| & \textrm{otherwise},
    \end{cases}
  \end{equation}
  and
  \begin{equation}\label{eq:umax2}
    \Delta_{max} = 
    \begin{cases}
      \epsilon_{1}\big|max\big(u_{i}, u_{i+1}\big) \big|  & \textrm{if } \sigma_{i-1} \sigma_{i+1} < 0      
                                                 \textrm{ and } \sigma_{i-1} > 0              
                                                \textrm{ or } \sigma_{i-1}\sigma_{i+1} \geq 0 
                                                 \textrm{ and } \sigma_{i-1} \sigma_{i} < 0 \\ 
      \epsilon_{0} \big|max\big(u_{i}, u_{i+1}\big) \big| & \textrm{otherwise}.
    \end{cases}
  \end{equation}
  The positive parameters $\epsilon_{0}$ and $\epsilon_{1}$, used for intervals with and without extrema, respectively, are introduced to adjust $\Delta_{min}$ and $\Delta_{max}$.
  This work extends the bounds in \cite{Ouermi2022} by introducing the parameter $\epsilon_{1}$ to allow for more flexibility on how to bound the interpolants in cases where an extremum is detected. 
  The choice for the positive parameters $\epsilon_{0}$ and $\epsilon_{1}$ depends on the underlying function and the input data used for the approximation. 
  As both $\epsilon_{0}$ and $\epsilon_{1}$ get smaller, the upper and lower bounds get tighter and the PPI method converges to the DBI method. 
  The choices for $\epsilon_{0}$ and $\epsilon_{1}$ are further discussed in Section \ref{subsec:software}.
  In Equation (\ref{eq:umin2}), the interval $I_{i}$ has a local maximum if $\sigma_{i-1} \sigma_{i+1}<0$ and $\sigma_{i-1} < 0$.
  Correspondingly, in Equation (\ref{eq:umax2}), the interval $I_{i}$ has a local minimum if $\sigma_{i-1} \sigma_{i+1}<0$ and $\sigma_{i-1} > 0$.
  In both Equations (\ref{eq:umin2}) and (\ref{eq:umax2}), the type of extremum is ambiguous if $\sigma_{i-1} \sigma_{i+1}$, and $\sigma_{i-1} \sigma_{i}<0$.
  Equation (\ref{eq:boundUP}) is equivalent to bounding $S_{n}(x)$ as follows:
  \begin{equation}\label{eq:boundSn}
    m_{\ell} \leq S_{n}(x) \leq m_r, 
  \end{equation} 
  where the factors $m_{\ell}$ and $m_{r}$ in Equation (\ref{eq:boundSn}) are expressed as 
  \begin{enumerate}%[label=\textbf{(\Roman*)}]
    \item : $u_{i+1} > u_{i}$
      \begin{equation}\label{eq:CaseI}
        m_{\ell} = min\bigg( 0, \frac{u_{min}-u_{i}}{u_{i+1}-u_{i}}\bigg), \textrm{ and }
        m_{r} = max\bigg( 1, \frac{u_{max}-u_{i}}{u_{i+1}-u_{i}}\bigg)
      \end{equation}
    \item : $u_{i+1} < u_{i}$
      \begin{equation}\label{eq:CaseII}
        m_{\ell} = min\bigg( 0, \frac{u_{max}-u_{i}}{u_{i+1}-u_{i}}\bigg), \textrm{ and }
        m_{r} = max\bigg( 1, \frac{u_{min}-u_{i}}{u_{i+1}-u_{i}}\bigg).
      \end{equation}
  \end{enumerate}
  The DBI method can be recovered from the PPI methods by setting $m_{\ell} = 0$ and $m_{r} = 1$.
  For $u_{i}=u_{i+1}$, $m_{\ell}$, $m_r$ and $U_{n}(x)$ as written in Equations (\ref{eq:CaseI}), (\ref{eq:CaseII}), and (\ref{eq:UN3}) are not defined.
    This limitation is addressed by re-writing $U_{n}(x)$ as 
    \begin{equation*}
      U_{n}(x) = u_{i} + U[x_{1}^{l}, \cdots, x_{1}^{r}](x_{i+1}-x_{i})(x_{1}^{r}-x_{1}^{l}) S_{n}(x),
    \end{equation*}
    where $S_{n}(x)$ is expressed as follows:
    \begin{equation*}
      S_{n}(x) = \sum_{j=1}^{n-1}\bar{s}_{j}.
    \end{equation*}
    The summation starts at $j=1$ because the linear term $\frac{u_{i+1}-u_{i}}{x_{i+1}-x_{i}}(x-x_{i}) =0$. 
    Let 
    \begin{equation*}
      w = U[x_{1}^{l}, \cdots, x_{1}^{r}](x_{i+1}-x_{i})(x_{1}^{r}-x_{1}^{l}).
    \end{equation*}
    $\bar{\lambda}_{j}$ in this context is defined as 
    \begin{equation*}
    \bar{\lambda}_{j} = \frac{U[x_{j}^{l}, \cdots, x_{j}^{r}]}{w}
    \prod_{k=0}^{j}(x_{k}^{r}-x_{k}^{l}).
  \end{equation*}
  For $u_{i}= u_{i+1}$, the parameters $m_{\ell}$ and $m_r$ in Equation (\ref{eq:boundSn}) are then defined according to
  \begin{enumerate}
    \item : $U[x_{1}^{l}, \cdots, x_{1}^{r}] >0$
      \begin{equation}\label{eq:ml_mr_1}
        m_{\ell} = min\bigg( 0, \frac{u_{min}-u_{i}}{w} \bigg), \textrm{ and }
        m_{r} = max\bigg( 1, \frac{u_{max}-u_{i}}{w}\bigg)
      \end{equation}
    \item : $U[x_{1}^{l}, \cdots, x_{1}^{r}] <0$
      \begin{equation}\label{eq:ml_mr_2}
        m_{\ell} = min\bigg( 0, \frac{u_{max}-u_{i}}{w}\bigg), \textrm{ and }
        m_{r} = max\bigg( 1, \frac{u_{min}-u_{i}}{w}\bigg).
      \end{equation}
  \end{enumerate}
    For $U[x_{i}, x_{i+1}] = U[x_{1}^{l}, \cdots, x_{1}^{r}] =0$, the data $u_{i-1}$, $u_{i}$, $u_{i+1}$, and $u_{i+2}$ have the same value ($u_{i-1}=u_{i}=u_{i+1}=u_{i+2}$). 
  In this case, the algorithm approximates the function in the interval $I_{i}$ with a linear interpolant.

  The positivity-preserving result in Equation (\ref{eq:boundUP}) is obtained by successively imposing bounds on the quadratic, cubic, and higher order terms in the expression of $S_{n}(x)$ defined in Equation (\ref{eq:SN1}).
  The reconstruction procedure begins by considering the linear and quadratic terms from $S_{n}(x)$ in Equation (\ref{eq:SN1}) and imposing the following bounds:
  \begin{equation}\label{eq:ppi_quad2}
    m_{\ell} \leq s + \frac{s(s-1)}{d_{1}}\bar{\lambda}_{1} \leq m_{r}.
  \end{equation}
  Equation (\ref{eq:ppi_quad2}) can be reorganized to obtain
  \begin{equation*}
    \bigg(\frac{m_{r}-1}{s(s-1)}-\frac{1}{s}\bigg)d_{1} 
    \leq \bar{\lambda}_{1} \leq 
    \bigg(\frac{m_{\ell}}{s(s-1)} - \frac{1}{(s-1)} \bigg)d_{1}.
  \end{equation*}
  Noting that $\frac{1}{s(s-1)} \leq -4$, $\frac{1}{s} \geq 1$, and $\frac{1}{s-1} \leq -1$, we obtain
  \begin{equation}\label{eq:ppi_lambda_bar_1}
    \bigg(-4(m_{r}-1)-1\bigg)d_{1} 
    \leq \bar{\lambda}_{1} \leq 
    \bigg(-4m_{\ell} + 1 \bigg)d_{1}.
  \end{equation}
  The bounds from Equation (\ref{eq:ppi_lambda_bar_1}) are extended to bound the cubic form by requiring that what multiplies $\bar{\lambda}_{1}$ must fit into the inequality in Equation (\ref{eq:ppi_lambda_bar_1}).  
  Thus, for the cubic case, Equation (\ref{eq:ppi_lambda_bar_1}) becomes 
  \begin{equation*} \label{eq:ppi_lambda_bar_2_2}
    \bigg(-4(m_{r}-1)-1\bigg)d_{1} \leq \bar{\lambda}_{1}\bigg(1 + \frac{(s-t_{2})}{d_{2}}\lambda_{2} \bigg)
                                   \leq d_{1}\bigg(-4m_{\ell} + 1 \bigg).
  \end{equation*}
  When $t_{2}$ defined in Equation (\ref{eq:tk}) is negative, $s-t_{2}$ has a maximum value at $s=1$ and a minimum value at $s=0$. 
  $\bar{\lambda}_{2}$ is then bounded by 
  \begin{equation*}\label{eq:ppi_lambda_bar_2}
    \frac{d_{2}}{(1-t_{2})}\bigg(\bigg(-4(m_{r}-1)-1\bigg)d_{1} -\bar{\lambda}_{1}\bigg)\leq\bar{\lambda}_{2}\leq \bigg(d_{1}\bigg(-4m_{\ell} + 1 \bigg)-\bar{\lambda}_{1}\bigg) \frac{d_{2}}{(1-t_{2})}.
  \end{equation*}
  When $t_{2}$ is positive, $\frac{1}{1-t_{2}}$ is substituted by $\frac{1}{-t_{2}}$ and the inequalities 
  $\leq$ with $\geq$ and vice versa are swapped. 

  This procedure is continued to quartic and higher order interpolants to produce the recursive expression for the bounds on $\bar{\lambda}_{j}$ for the PPI and DBI methods as follows:
  \begin{subequations}
   \begin{equation}\label{eq:B_PPIminus}
        B_{j}^{-} = 
      \begin{cases}
        (-4(m_{r}-1) -1)d_{1} \quad j =1 \\
        (B_{j-1}^{-} -\bar{\lambda}_{j-1}) \frac{d_{j}}{1-t_{j}}, \textrm{ if } t_{j} \in (-\infty, 0] \quad j>1 \\ 
        (B_{j-1}^{+} -\bar{\lambda}_{j-1}) \frac{d_{j}}{-t_{j}}, \textrm{ if } t_{j} \in (0, +\infty) \quad j>1, 
      \end{cases}
    \end{equation}
    \textrm{and}
    \begin{equation}\label{eq:B_PPIplus}
      B_{j}^{+} =
      \begin{cases}
         (-4m_{\ell}+1)d_{1} , \quad j = 1\\
         (B_{j-1}^{+} -\bar{\lambda}_{j-1}) \frac{d_{j}}{1-t_{j}}, \textrm{ if } t_{j} \in (-\infty, 0] \quad j >1\\ 
         (B_{j-1}^{-} -\bar{\lambda}_{j-1}) \frac{d_{j}}{-t_{j}}, \textrm{ if } t_{j} \in (0, +\infty) \quad j >1.
      \end{cases}
    \end{equation}
  \end{subequations}
  $B_{1}^{-}$  and $B_{1}^{+}$ are defined as $-d_{1}$ and $d_{1}$ for the DBI method, whereas for the PPI method, they are defined as $(-4(m_{r}-1) -1)d_{1}$ and $(-4m_{\ell}+1)d_{1}$, respectively.
  We refer the reader to Theorems 1 and 2 in ~\cite{Ouermi2022} for more details on the mathematical foundation used to build the positivity-preserving software.

  \section{Algorithms and Software}
\label{sec:algorithm}
  This section describes the algorithms and different components used in the data-bounded and positivity-preserving software.
  The software developed in this work provides 1D, 2D, and 3D implementations of the DBI and PPI methods for uniform and nonuniform structured meshes.
  The 1D implementation is constructed based on the mathematical framework provided in Section \ref{sec:background}.
  The 2D and 3D implementations are obtained via a tensor product of the 1D version.

  \subsection{Algorithms}
  \label{subsec:algorithms}
    The algorithms provide the necessary elements to construct the data-bounded or positive interpolants.
    Rogerson and Meiburg ~\cite{rogerson1990numerical} showed that the ENO reconstruction can lead to a left- or right-biased stencil that causes stability issues when used to solve hyperbolic equations. 
    Shu ~\cite{shu1990numerical} addressed this limitation by introducing a bias coefficient used to target a preferred stencil. 
    As indicated in ~\cite{Ouermi2022}, the left- and right-biased stencil can fail to recover hidden extrema. For a given interval $I_{i}$,  the left- and right-biased stencil does not include the points $x_{i-1}$ or $x_{i+1}$, respectively.
    \textbf{Algorithm I} addresses these limitations by extending the algorithm in ~\cite{Ouermi2022} to introduce more options for the adaptive stencil selection process described below.
    In addition to the symmetry-based points selection in ~\cite{Ouermi2022}, \textbf{Algorithm I} includes ENO-type and locality-based point selection processes.

    At any given step $j$, the next point inserted into $\mathcal{V}_{j}$ can be to the left or right.
    Let $x_{p}$ and $x_{q}$ be the mesh points immediately to the left and right of $\mathcal{V}_{j}$, respectively.
    
    We define $\bar{\lambda}_{j+1}^{-}$ and $\bar{\lambda}_{j+1}^{+}$ as follows:   
    \begin{equation}\label{eq:lambda_plus_minus}
      \begin{cases}
        \bar{\lambda}_{j+1}^{-} =\bar{\lambda}_{j+1} & \textrm{ with } \mathcal{V}_{j+1} = \{x_{p} \} \cup \mathcal{V}_{j} \\ 
        \bar{\lambda}_{j+1}^{+} =\bar{\lambda}_{j+1} & \textrm{ with } \mathcal{V}_{j+1} = \mathcal{V}_{j} \cup \{x_{q} \}. 
      \end{cases}
    \end{equation}
    
    The terms $\bar{\lambda}_{j+1}^{-}$ and $\bar{\lambda}_{j+1}^{+}$ in Equation \ref{eq:lambda_plus_minus} correspond to the case where the stencil inserted is to the left and right, respectively.
    Given $\mathcal{V}_{j}$, let $\mu_{j}^{l}$ be the number of points to the left of $x_{i}$ and $\mu_{j}^{r}$ the number of points to the right.
    \textbf{Algorithm I} extends the algorithm in ~\cite{Ouermi2022} by introducing a user-supplied parameter $st$ used to guide the procedure for stencil construction.
    In the cases where adding both $x_{p}$ (to the left) or $x_{q}$ (to the right) are valid, the algorithm makes the selection based on the three cases below:
    \begin{itemize}
      \item If $st=1$ (default), the algorithm chooses the point with the smallest divided difference, as in the ENO stencil. 
      \item If $st=2$, the point to the left of the current stencil is selected if the number of points to the left of $x_{i}$ is smaller than the number of points to the right. 
      Similarly, the point to the right is selected if the number of points to the right of $x_{i}$ is smaller than the number of points to the left.
      When both the number of points to the right and left are the same, the algorithm chooses the point with the smallest $\bar{\lambda}_{j+1}$.
      \item If $st=3$, the algorithm chooses the point that is closest to the starting interval $I_{i}$.
      It is important to prioritize the closest points in cases where the intervals surrounding $I_{i}$ vary significantly in size.
      These variations are found in computational problems for which different resolutions are used for different parts of the domain.
    \end{itemize}
    \textbf{Algorithm II} describes the 1D DBI and PPI methods built using the mathematical framework in Section \ref{sec:background} and \textbf{Algorithm I}.
    \textbf{Algorithm II} further extends the constraints in \cite{Ouermi2022} by introducing the user-supplied positive parameter $\epsilon_{1}$ that is used to impose upper and lower bounds on the interpolants according to Equations (\ref{eq:umin2}) and (\ref{eq:umax2}).
    The positive parameters $\epsilon_{0}$ and $\epsilon_{1}$ are used for intervals without and with an extremum, respectively.
    The user-supplied parameter $im$ is used to choose between the DBI and PPI methods.
    \textbf{Algorithm III} and \textbf{Algorithm IV} describe the extension from 1D to 2D and 3D, respectively.
    Both \textbf{Algorithm III} and \textbf{IV} are constructed by successively applying \textbf{Algorithm II} to each dimension.
    Given that the DBI and PPI methods are nonlinear, the order in which \textbf{Algorithm II} is used can lead to different approximation results.
    In this paper and in \cite{Ouermi2022}, the 1D DBI and PPI are first applied to $x$, then the $y$ dimension, and finally the  $z$ dimensions, as indicated in \textbf{Algorithm III} and \textbf{IV}.
    In \textbf{Algorithm III} and \textbf{IV}, the input and intermediate data values are modified only by \textbf{Algorithm I}, which preserves data-boundedness or positivity.
    Therefore, the resulting solutions from the 2D and 3D extensions will preserve data-boundedness and positivity.
    Similar to \textbf{Algorithm II}, the choices of parameters $st$, $\epsilon_{0}$, and $\epsilon_{1}$ influence the quality of the approximation in \textbf{Algorithm III} and \textbf{IV}.
    In the 1D, 2D, and 3D cases, the choice for parameters $st$, $\epsilon_{0}$, and $\epsilon_{1}$ is dependent on the data.
    In the case of 2D and 3D, applying the 1D PPI method (\textbf{Algorithm II}) along the $x$ and/or $y$ dimensions may introduce oscillations that could be amplified when applying the 1D PPI in the subsequent dimensions.
    These oscillations can be significantly reduced with small values for $\epsilon_{0}$ and $\epsilon_{1}$.
    The parameters $\epsilon_{0}$ and $\epsilon_{1}$ should be chosen to be small enough such that hidden extrema are recovered without introducing new large oscillations.
    In \textbf{Algorithm II} and \textbf{IV}, $\mathcal{M}_{\square x \times \square y}$ and $\mathcal{M}_{\square x \times \square y \square z}$ represent 2D and 3D meshes obtained by taking the tensor product of the 1D mesh along the $x$ and $y$ dimensions for the 2D mesh, and along the $x$, $y$, and $z$ dimensions for the 3D mesh. 
    The square $\square$ represents $n$ or $m$.
    
    \textbf{\underline{Algorithm I}} \\
    \textbf{Input:} $\mu_{j}^{l}$, $\mu_{j}^{l}$, $x_{p}$, $x_{i}$, $x_{q}$, $x_{i+1}$, $U[x_{p}, \cdots, x_{j}^{r}]$, $U[x_{j}^{l}, \cdots, x_{q}]$ $\bar{\lambda}_{j+1}^{-}$, $\bar{\lambda}_{j+1}^{+}$, and $st$.
    \begin{enumerate}
      \item if $st = 1$ 
        \begin{itemize}
           \item if $|U[x_{p} \cdots, x_{j}^{r}]| < |U[x_{j}^{l}, \cdots, x_{q}]|$, then insert a new stencil point to the left;  
           \item else if $|U[x_{p} \cdots, x_{j}^{r}]| > |U[x_{j}^{l}, \cdots, x_{q}]|$, then insert a new stencil point to the right; 
           \item else  insert a new stencil point to the right if $|\bar{\lambda}_{j+1}^{-}|\geq|\bar{\lambda}_{j+1}^{+}|$; otherwise, insert a new point to left;
        \end{itemize}
      \item if $st = 2$ 
        \begin{itemize}
           \item if $\mu_{j}^{l} < \mu_{j}^{r}$, then insert a new stencil point to the left;  
           \item else if $\mu_{j}^{l} > \mu_{j}^{r}$, then insert a new stencil point to the right; 
           \item else  insert a new stencil point to the right if $|\bar{\lambda}_{j+1}^{-}|\geq|\bar{\lambda}_{j+1}^{+}|$; otherwise, insert a new point to left;
        \end{itemize}
      \item else $st = 3$
        \begin{itemize}
           \item if $|x_{p}-x_{i}| < |x_{q}-x_{i+1}|$, then insert a new stencil point to the left;  
           \item else if $|x_{p}-x_{i}| > |x_{q}-x_{i+1}|$, then insert a new stencil point to the right; 
           \item else  insert a new stencil point to the right if $|\bar{\lambda}_{j+1}^{-}|\geq|\bar{\lambda}_{j+1}^{+}|$; otherwise, insert a new point to left;
        \end{itemize}
    \end{enumerate} 
    \textbf{\underline{Algorithm II (1D)}} \\
    \textbf{Input:} $\{x_{i}\}_{i=0}^{n}$, $\{u_{i}\}_{i=0}^{n}$, $\{\tilde{x}_{i}\}_{i=0}^{\tilde{n}}$, $d$, $st$ $\epsilon_{0}$, $im$, and $\epsilon_{1}$.
    \textbf{Output:} $\{\tilde{u}_{i}\}_{i=0}^{\tilde{n}}$.
    \begin{enumerate}
       \item Select an interval $[x_{i}, x_{i+1}]$.  Let $\mathcal{V}_{0}=\{x_{i}, x_{i+1}\}=\{x_{0}^{l}, x_{0}^{r}\}$.
       \item If $\sigma_{i-1}\sigma_{i+1} < 0$ or $\sigma_{i-1}\sigma_{i} < 0$, then the interval $I_{i}$ has a hidden local extremum.
             For the boundary intervals, we assume that the divided differences to the left and right have the same sign.
       \item Compute $u_{min}$ and $u_{max}$ using Equations (\ref{eq:uminumax}), (\ref{eq:umin2}), and (\ref{eq:umax2}).
       \item Compute $m_{r}$ and $m_{\ell}$ based on Equations (\ref{eq:CaseI}) and (\ref{eq:CaseII}) or Equations (\ref{eq:ml_mr_1}) and (\ref{eq:ml_mr_2}). 
             For DBI, set $m_{r}=1$ and $m_{\ell} = 0$.
         \item Given a stencil $\mathcal{V}_{j}$, 
         \begin{itemize}
           \item if $B^{-}_{j+1} \leq \bar{\lambda}_{j+1}^{+} \leq B_{j+1}^{+}$ and $B^{-}_{j+1} \leq \bar{\lambda}_{j+1}^{-} \leq B_{j+1}^{+}$, choose the point to add based on \textbf{Algorithm I}
          \item else if $B^{-}_{j+1} \leq \bar{\lambda}_{j+1}^{-} \leq B_{j+1}^{+}$, then insert a new stencil point to the left; 
           \item else if $B^{-}_{j+1} \leq \bar{\lambda}_{j+1}^{+} \leq B_{j+1}^{+}$, then insert a new stencil point to the right; 
         \end{itemize}
       \item This process (Steps $3$) iterates until the halting criterion that the ratio of divided differences lies outside the required bounds stated above or the stencil has $d+1$ points, with $d$ being the target degree for the interpolant. 
       \item Evaluate the final interpolant $U^{l}(x)$ (for DBI) or $U^{p}(x)$ (for PPI) at the output points $\tilde{x}_{i}$ that are in $I_{i}$.
       \item Repeat Steps $1$--$7$ for each interval in the input 1D mesh.
     \end{enumerate} 
    \textbf{\underline{Algorithm III (2D)}}\\
    \textbf{Input:}$\{x_{i}\}_{i=0}^{nx}$,$\{y_{i}\}_{i=0}^{ny}$ $\{u_{i,j}\}_{i,j=0}^{nx, ny}$, $d$, $st$ $\epsilon_{0}$, $\{\tilde{x}_{i}\}_{i=0}^{mx}$, $\{\tilde{y}_{j}\}_{j=0}^{my}$, $im$, and $\epsilon_{1}$.
    \textbf{Output:} $\{\tilde{u}_{i,j}\}_{i,j=0}^{mx, my}$.
    \begin{enumerate}
        \item  Interpolate from the mesh $\mathcal{M}_{nx \times ny}$ to $\mathcal{M}_{mx \times ny}$ by applying the \textbf{Algorithm II (1D)} along the x dimension for each index $j$ along the y direction. 
        More precisely, for each index $j$ ($0 \leq j \leq nx$), \textbf{Algorithm II (1D)} uses $\{x_{i}\}_{i=0}^{nx}$, $\{u_{i,j}\}_{i=0}^{nx}$ to interpolate from $\{x_{i}, y_{j}\}_{i=0}^{nx}$ to $\{\tilde{x}_{i}, y_{j}\}_{i=0}^{mx}$. 
        The interpolation solution is saved in $\{q_{i,j}\}_{i,j=0}^{mx, ny}$
        \item Use the solution$\{q_{i,j}\}_{i,j=0}^{mx, ny}$ from step (1) to interpolate from $\mathcal{M}_{mx \times ny}$ to the final output mesh $\mathcal{M}_{mx \times my}$ by applying the \textbf{Algorithm II (1D)} along the y dimension  for each index $i$ along the x direction. 
        For each index $i$ ($0<i<mx$), \textbf{Algorithm II (1D)} uses ($\{y_{j}\}_{j=0}^{ny}$, $\{q_{i,j}\}_{j=0}^{ny}$) to interpolate from $\{\tilde{x}_{i}, y_{j}\}_{j=0}^{ny}$ to $\{\tilde{x}_{i}, \tilde{y}_{j}\}_{j=0}^{my}$. 
        The interpolation final results are saved in $\{\tilde{u}_{i,j} \}_{i,j}^{mx,my}$
    \end{enumerate}
    \textbf{\underline{Algorithm IV (3D)}}\\
    \textbf{Input:}$\{x_{i}\}_{i=0}^{nx}$,$\{y_{i}\}_{i=0}^{ny}$, $\{z_{k}\}_{k=0}^{nz}$, $\{u_{i,j,k}\}_{i,j,k=0}^{nx,ny,nz}$, $d$, $st$ $\epsilon_{0}$, $\{\tilde{x}_{i}\}_{i=0}^{mx}$, $\{\tilde{y}_{j}\}_{j=0}^{my}$, $\{\tilde{z}_{k}\}_{k=0}^{mz}$, $im$, and $\epsilon_{1}$.
    \textbf{Output:} $\{\tilde{u}_{i,j,k}\}_{i,j=0}^{mx,my,mz}$.
    \begin{enumerate}
        \item  Interpolate from the mesh $\mathcal{M}_{nx \times ny \times nz}$ to $\mathcal{M}_{mx \times ny \times nz}$ by applying the  \textbf{Algorithm II (1D)} along the x dimension  for each pair $j,k$. 
        More precisely, for each each pair $j,K$ ($0 \leq j \leq nx$ and $0 \leq k \leq nz$), \textbf{Algorithm II (1D)} uses $\{x_{i}\}_{i=0}^{nx}$, $\{u_{i,j,k}\}_{i=0}^{nx}$ to interpolate from $\{x_{i}, y_{j}, z_{k}\}_{i=0}^{nx}$ to $\{\tilde{x}_{i}, y_{j}, z_{k}\}_{i=0}^{mx}$. 
        The interpolation solution is saved in $\{q_{i,j,k}\}_{i,j,k=0}^{mx, ny, nz}$
        \item Use the solution$\{q_{i,j,k}\}_{i,j, k=0}^{mx, ny, nz}$ from step (1) to interpolate from $\mathcal{M}_{mx \times ny \times nz}$ to $\mathcal{M}_{mx \times my \times nz}$ by applying the \textbf{Algorithm II (1D)} along the $y$ dimension. 
        For each pair $i,k$ ($0\leq i\leq mx$ and $0 \leq k \leq nz$), \textbf{Algorithm II (1D)} uses ($\{y_{j}\}_{j=0}^{ny}$, $\{q_{i,j,k}\}_{j=0}^{ny}$) to interpolate from $\{\tilde{x}_{i}, y_{j}, z_{k}\}_{j=0}^{ny}$ to $\{\tilde{x}_{i}, \tilde{y}_{j}, z_{k}\}_{j=0}^{my}$. 
        The interpolation results are saved in $\{\tilde{g}_{i,j,k} \}_{i,j}^{mx,my,nz}$.
        \item Use the solution$\{q_{i,j,k}\}_{i,j, k=0}^{mx, ny, nz}$ from step (1) to interpolate from $\mathcal{M}_{mx \times ny \times nz}$ to $\mathcal{M}_{mx \times my \times nz}$ by applying the \textbf{Algorithm II (1D)} along the $y$ dimension. 
        For each pair $i,j$ ($0\leq i\leq mx$ and $0 \leq j \leq my$), \textbf{Algorithm II (1D)} uses ($\{z_{k}\}_{k=0}^{nz}$, $\{g_{i,j,k}\}_{k=0}^{nz}$) to interpolate from $\{\tilde{x}_{i}, \tilde{y}_{j}, z_{k}\}_{k=0}^{nz}$ to $\{\tilde{x}_{i}, \tilde{y}_{j}, \tilde{z}_{k}\}_{k=0}^{mz}$. 
        The interpolation final results are saved in $\{\tilde{u}_{i,j,k} \}_{i,j}^{mx,my,mz}$.
    \end{enumerate}

    \subsection{Software Description}
    \label{subsec:software}
    The DBI and PPI software implementation is guided by the algorithms described above. 
    HiPPIS is available at \url{https://github.com/ouermijudicael/HiPPIS}.
    The software can be organized into four major parts: (1) computation of divided differences, (2) calculations of upper and lower bounds for each interval, (3) a stencil construction procedure, and (4) 1D, 2D, and 3D DBI and PPI implementations.

    The divided differences are essential to the DBI and PPI methods because they are used in the calculations of $\bar{\lambda}_{j}$ and the stencil selection process.
    The divided differences are computed using the standard recurrence form in Equation (\ref{eq:divDiff}) and stored in a table of dimension $n\times(d+1)$ where $d$ is the maximum polynomial degree for each interpolant.
    Given that the maximum degree is $d$, it is sufficient to consider the $d+1$ divided differences for the stencil selection process and the construction of the final polynomial interpolant for each interval.

    The bounds on each interpolant are obtained from Equation (\ref{eq:uminumax}), (\ref{eq:umin2}), and (\ref{eq:umax2}) where the positive parameters $\epsilon_{0}$  and $\epsilon_{1}$ are user-supplied values used to adjust the bounds for the interval with and without extremum, respectively.
    The adjustment focuses on removing large oscillations as much as possible while still allowing high-degree polynomial interpolants that meet the positivity requirements.  
    The stencil selection process requires the computation of $B_{j}^{+}$ and $B_{j}^{-}$, which are both dependent on $d_{j}$, $t_{j}$, and $\bar{\lambda}_{j}$.
    The stencil $\mathcal{V}_{j}$ is constructed from $\mathcal{V}_{j-1}$ by appending a point to the left or right of $\mathcal{V}_{j-1}$.
    When appending to either the right or left meets the requirements for positivity, the software offers three possible options for choosing from both points that can be set by the user.
    The first and default option ($st=1$) chooses the stencil with the smallest divided difference, similar to the ENO-like approach.
    The second option ($st=2$) prioritizes the choice that makes the stencil more symmetric around $x_{i}$.
    The third option ($st=3$) chooses the point closest to the starting interval $I_{i}$, thus prioritizing locality.

    The 1D DBI and PPI methods use (1), (2), and (3) as building blocks to construct the final approximation.
    Once the final stencil has been selected, the interpolant is built using a Newton polynomial representation and then evaluated at the corresponding output points.
    The Newton polynomial is used here because its coefficients/divided differences are available.
    %
    %Given that the coefficients are available, constructing the final interpolant using a Newton polynomial is computationally less expensive than using a Lagrange interpolation.
    %
    The 2D and 3D implementations successively use the 1D version along each dimension to construct the final approximation on uniform and nonuniform structured meshes.

    The interfaces for the 1D, 2D, and 3D DBI and PPI subroutines are designed to be similar to widely used interfaces for polynomial interpolation such as PCHIP, and can be incorporated into larger application codes.
    The interfaces require 
    \begin{itemize} 
      \item the input mesh points and the data values associated with those points, 
      \item the maximum polynomial degree to be used for each interpolant,
      \item the interpolation method to be used (DBI or PPI), and  
      \item the output mesh points.
    \end{itemize}
    Listing \ref{code:interface} shows examples of how to use the 1D, 2D, and 3D interfaces for DBI and PPI in HiPPIS. 
    The variables \textit{x, y,} and \textit{z} are 1D arrays used to define the input meshes, and \textit{xout, yout,} and \textit{zout} are used to define the output meshes. 
    The variables \textit{v, v2D,} and \textit{v3D} correspond to the input data values associated with the input meshes. 
   The parameters \textit{d} and \textit{im} (1, or 2) are used to indicate the target polynomial degree and the interpolation method to be used. 
    For DBI and PPI, the parameter \textit{im} is set to 1 and 2, respectively. 
    The parameters $st$, $\epsilon_{0}$, and $\epsilon_{1}$ are optional parameters that are set to $3$, $0.01$, and $1$ by default, as explained below.
    The choice of the optional parameters depends on the underlying function and the input data.  

    In problems for which different resolutions are used for different parts of the computational domain, \textit{st=3} is a preferable choice. 
    The algorithm prioritizes the closest points to the starting interval $I_{i}$ if \textit{st=3}. 
    This choice is particularly important in regions where the size of the intervals varies significantly. 
    For cases when smoothness is the primary goal, \textit{st=1} is a suitable choice. 
    For \textit{st=1}, the \textbf{Algorithm I} prioritizes smoothness by choosing the points with the smallest divided differences during the stencil construction process. 
    Both the \textit{st=1} and \textit{st=3} can lead to a left- or right-biased stencil. 
    In these instances, \textit{st=2} can be used to remove the bias. 
    For \textit{st=2}, the algorithm prioritizes a symmetric stencil. 
    The default value of \textit{st} is set to 3 because the examples in this study indicate that st=3 leads to better approximations compared to st=1 or 2 and locality is often a highly desired property in many computational problems. 

    The positive parameters $\epsilon_{0}$ and $\epsilon_{1}$ are used to bound the interpolants for the intervals with and without extrema, respectively. 
    The configurations in ~\cite{Ouermi2022} and ~\cite{arkivtajo} correspond to setting the parameters $\epsilon_{0}$  and $\epsilon_{1}$ to the default values of $0.01$ and $1$, respectively. 
    The values of $\epsilon_{0}$ and $\epsilon_{1}$ are chosen such that the lower and upper bounds on each interpolant are relaxed enough to allow for a high-order polynomial that does not introduce undesirable oscillations. 
    For profiles that are prone to oscillation such as the logistic functions, it is important to choose small values for $\epsilon_{0}$ and $\epsilon_{1}$. 
    For $N \times N = 17 \times 17$, the approximation leads to large oscillations if $\epsilon_{0}$ and $\epsilon_{1}$ are greater than $10^{-4}$. 
    For intervals without extrema, it is important to keep $\epsilon_{0}$ small to not introduce new extrema. 
    For the intervals with extrema, $\epsilon_{1}$ needs to be large enough to allow for recovery of hidden extrema but small enough to not cause undesired large oscillations. 
    This is very challenging given that the sizes of the peaks are not known a priori. 
    The default values of $\epsilon_{1}=1$ are such that the interpolant maximum value is twice $max(u_{i}, u_{I+1})$. 
    This default value of one is sufficient for the modified Runge and TWP-ICE examples. 
    However, in the case of BOMEX, smaller values of $\epsilon_{1}\leq10^{-5}$ are required to remove undesired oscillations. 
    In practice, it is prudent to start with a small value for $\epsilon_{0}$ and $\epsilon_{1}$ and increase them as needed if the approximation fails to recover hidden extrema or uses low-degree polynomial interpolants. 

\begin{figure}[H]
     \centering
     \includegraphics[scale = 0.27]{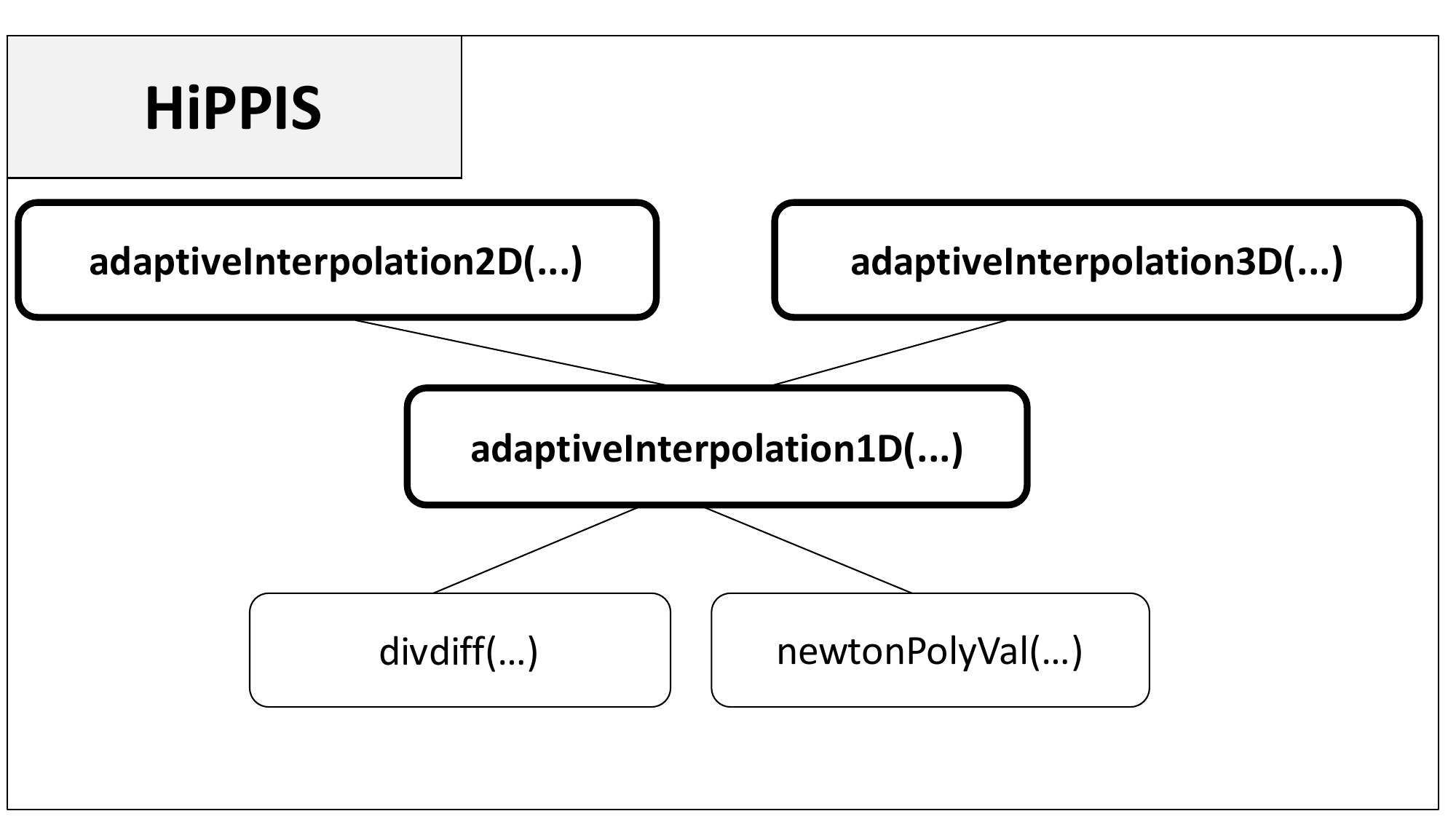}
     \caption{Diagram showing the components of the main module used to build the HiPPIS software.}
     \label{fig:HPPIS-diagram}
   \end{figure}

    Figure \ref{fig:HPPIS-diagram} is a diagram of the different components of the main module of HiPPIS.
    The function \textit{divdiff(...)} is used to calculate the divided differences needed for \textbf{Algorithms I} and \textbf{II}.
    Once the final stencil is constructed, the function \textit{newtonPolyVal(...)} is used to build and evaluate the positive interpolant at the corresponding output points.
    The major part of the data-boundedness and positivity preservation including \textbf{Algorithms I} and \textbf{II} is in the function \textit{adaptiveInterpolation1D(...)}.
    This function is used for the 1D approximation or mapping problems and depends on the function \textit{divdiff(...)} and \textit{newtonPolyVal(...)}.
    The functions \textit{adaptiveIterpolation2D(...)} and \textit{adaptiveInterpolation3D(...)} use \textit{adaptiveInterpolation1D(..)} to construct the data-bounded or positive polynomial approximations on 2D and 3D structured tensor product meshes, respectively.
    The interfaces for the 1D, 2D, and 3D interpolations, in bold, require the parameter $im$, which is used to indicate the interpolation method chosen.
    For the DBI and PPI methods, the parameter $im$ is set using 1 and 2, respectively. 
    HiPPIS does not allow for any other choices for the parameter $im$. 
      \begin{lstlisting}[language=Matlab, mathescape=true, caption=Interface examples, label=code:interface]
% 1D example 
vout = adaptiveInterpolation1D(x, v, xout, d, im, st, $\epsilon_{0}$, $\epsilon_{1}$ );

%2D example
vout2D = adaptiveInterpolation2D(x, y, v2D, xout,yout, d, im, st, $\epsilon_{0}$, $\epsilon_{1}$ );

%3D example
vout3D = adaptiveInterpolation3D(x, y, z, v3D, xout, yout, zout, d, im, st, $\epsilon_{0}$, $\epsilon_{1}$ );
   \end{lstlisting}

    \section{Numerical Examples for Function Approximation}
  \label{sec:numerical-examples}
    This section provides 1D and 2D numerical examples used to evaluate the PCHIP, MQSI, DBI, and PPI methods.
    The results are based on the Fortran implementation of the software.
    These examples include a subset of the full suite of test problems considered in ~\cite{arkivtajo}.
    The interpolation methods are used to approximate positive functions from provided data values that are obtained by evaluating the 1D and 2D functions on a given set of mesh points.
    The function approximations in 2D examples use the 2D extension of the DBI and PPI methods described in \textbf{Algorithm III}.
    The 2D PCHIP and MQSI methods are obtained by successively applying the 1D PCHIP and MQSI algorithms to each dimension, similar to the steps described in \textbf{Algorithm III}.
    The 1D PCHIP and MQSI are first applied to $x$, then the $y$ dimension, and finally the $z$ dimension.
    Similar to 2D and 3D DBI and PPI, the order in which 1D PCHIP or MQSI is applied can lead to different results, as both methods are nonlinear and non-commutative. 
    Using a standard polynomial interpolation to approximate the different functions leads to negative values and oscillations. 
    In this section, the $L^{2}$-norms in the tables below are approximated using the trapezoidal rule with, $10^{4}$ and $10^{3} \times 10^{3}$ uniformly spaced points for the 1D and 2D examples, respectively. 
    For the numerical examples in Sections \ref{subsec:f1} - \ref{subsec:twp-ice}, the errors from using $st=1, 2$, and $3$ are similar, with $st=3$ leading to slightly smaller errors compared to $st=1$ and $st=2$. 
    Given that the results are similar, Tables \ref{tab:f1} - \ref{tab:twp-ice-map} show errors with the parameter $st$ set to $3$. 
    For the BOMEX example, the errors from the three choices are significantly different. 
    Therefore, the results from all three choices are included. 
    More test examples can be found in ~\cite{arkivtajo}.

    \subsection{Example I: Modified Runge Function}
    \label{subsec:f1}
      This example uses a modified version of the canonical Runge function defined as 
      \begin{equation}\label{eq:runge-1d}
         f_{1}(x) = \frac{0.1}{0.1+ 25x^{2}}, \quad x \in [-1,1].
      \end{equation}
      Approximating the modified Runge function in Equation (\ref{eq:runge-1d}) with a global standard polynomial leads to large oscillations.
      Table \ref{tab:f1} shows the $L^{2}$-errors norms when using the PCHIP, MQSI, DBI, and PPI methods to approximate $f_{1}(x)$. 
      The DBI and PPI methods lead to better approximation results compared to the PCHIP and MQSI  methods. 
      The errors from PCHIP and MQSI are comparable.
      In this example, the algorithms used in MQSI to approximate and adjust the derivatives to enforce monotonicity do not produce more accurate results compared to PCHIP.
      As the target polynomial degree increases from $d=4$ to $d=8$, the DBI approximation does not improve significantly compared to the PPI method.
      The relaxed nature of the PPI method allows for higher degree polynomial interpolants compared to DBI, PCHIP, and MQSI which lead to better approximations.

      \begin{table}[H]
        \centering
        \begin{tabular}{ c c c c c c c c c c c c c}
          \hline
          \hline
          $N$   && PCHIP    && MQSI &&\multicolumn{3}{c}{DBI} && \multicolumn{3}{c}{PPI}   \\

                &&$\mathcal{P}_{3}$ && $\mathcal{P}_{5}$ && $\mathcal{P}_{3}$ & $\mathcal{P}_{4}$ & $\mathcal{P}_{8}$
                                    &&$\mathcal{P}_{3}$ & $\mathcal{P}_{4}$ & $\mathcal{P}_{8}$ \\
          \hline
          17       && 3.99E-2  && 3.63E-2 &&  5.10E-2  &  2.91E-2  &  4.61E-2  &&  5.10E-2  &  2.91E-2  &  4.61E-2   \\
          33       && 4.52E-3  && 4.32E-3 &&  6.31E-3  &  9.57E-3  &  3.05E-3  &&  6.31E-3  &  9.57E-3  &  3.05E-3   \\
          65       && 2.79E-3  && 2.67E-3 &&  2.44E-3  &  2.49E-3  &  1.33E-3  &&  2.44E-3  &  2.49E-3  &  9.92E-4   \\
          129      && 6.23E-4  && 6.71E-4 &&  2.22E-4  &  1.21E-4  &  1.05E-4  &&  2.22E-4  &  1.21E-4  &  2.43E-5   \\
          257      && 1.17E-4  && 9.89E-5 &&  1.51E-5  &  1.15E-5  &  1.07E-5  &&  1.51E-5  &  4.68E-6  &  9.89E-8   \\
          \hline
          \hline
        \end{tabular}                                                                                  
        \caption{$L^2$-errors when using the PCHIP, MQSI, DBI, and PPI  methods to approximate the function $f_{1}(x)$.
                 $N$ represents the number of input points used to build the approximation.
                 The parameters $\epsilon_{0}$, $\epsilon_{1}$, and $st$ are set to $0.01, 1.0,$ and $3$, respectively.}
        \label{tab:f1}
      \end{table}
     
    \subsection{Example II: 1D Logistic Function}
      This test case uses a logistic function defined as 
      \begin{equation}\label{eq:smoothed-heaviside-1d}
        f_{2}(x) = \frac{1}{1+e^{-2kx}}, \quad k=100, \textrm{ and } x \in [-0.2, 0.2].
      \end{equation}
      This function is a smoothed analytical approximation of the Heaviside step function. 
      The logistic function in Equation (\ref{eq:smoothed-heaviside-1d}) is challenging because of the steep gradient at about $x=0$.
      Approximating $f_{2}(x)$ with a standard polynomial interpolation leads to large oscillations to the left and right of the gradient.
      In addition, the oscillations to the left produce negative values.

      Table \ref{tab:f2} shows $L^{2}$-error norms when using the PCHIP, MQSI, DBI, and PPI methods to approximate the smoothed logistic function $f_{2}(x)$.
      The MQSI method has larger errors compared to the other methods.
      In this case, the algorithms employed by MQSI to approximate and adjust derivatives values used to construct the monotonic quintic splines are less accurate than the one used in PCHIP.
      For a target polynomial degree $d=3$, the approximation errors using PCHIP, DBI, and PPI are comparable.
      Increasing the target polynomial degree improves the approximations for DBI and PPI, as shown in Table \ref{tab:f2}.
      The errors from both the DBI and PPI methods are similar because the logistic example has no hidden extrema, and the stencils used for both methods are the same, around $x=0$.
      The global error is dominated by the local errors in the region with steep gradients around $x=0$.

\begin{figure}
        \centering
        \begin{subfigure}{0.9 \textwidth}
            \centering
            \includegraphics[scale=0.35]{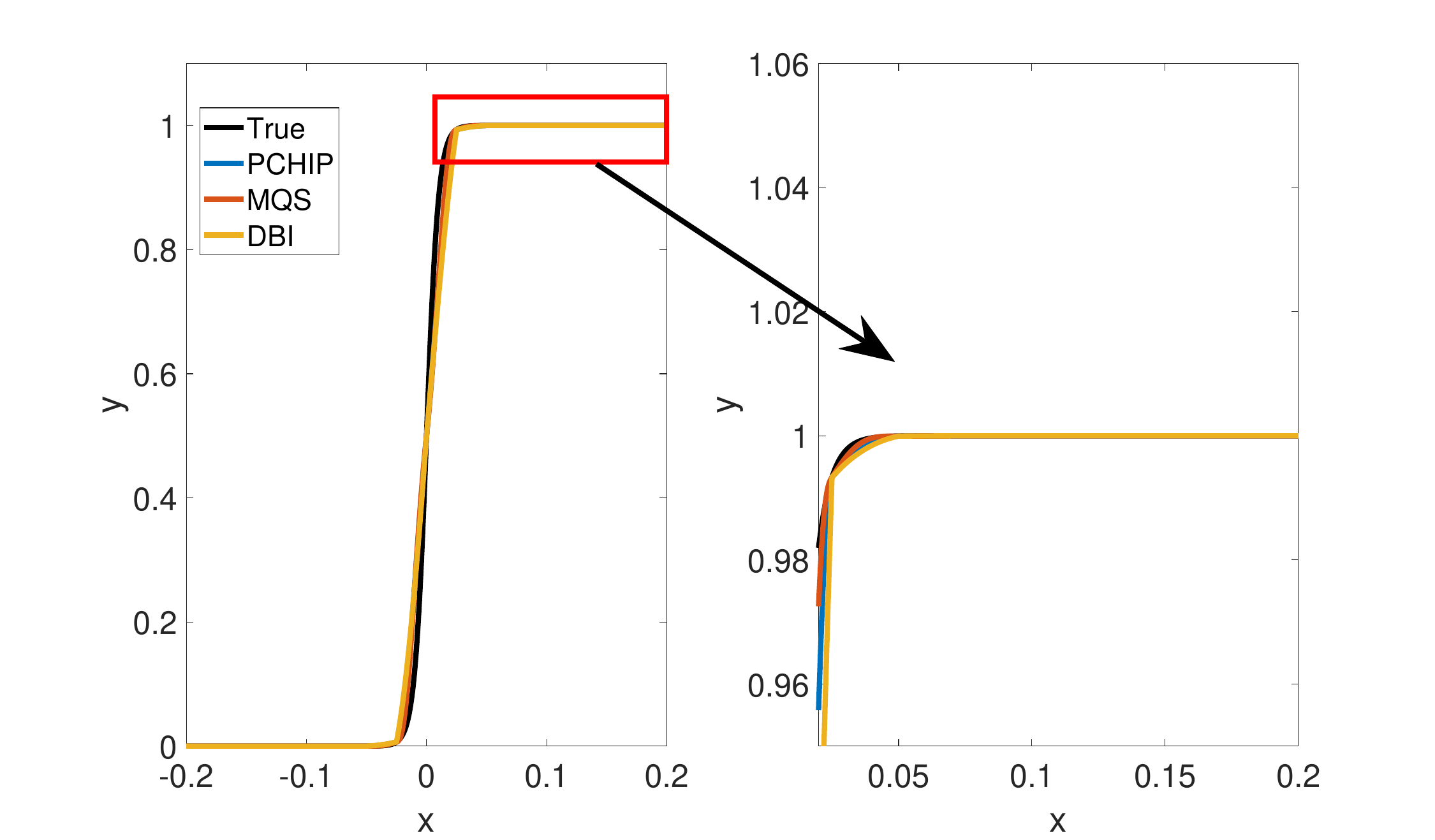}
            \subcaption{PCHIP, MQSI, and DBI with $d=8$.}
            \label{subfig:heaviside1}
        \end{subfigure}
        \begin{subfigure}{0.9 \textwidth}
            \centering
            \includegraphics[scale=0.35]{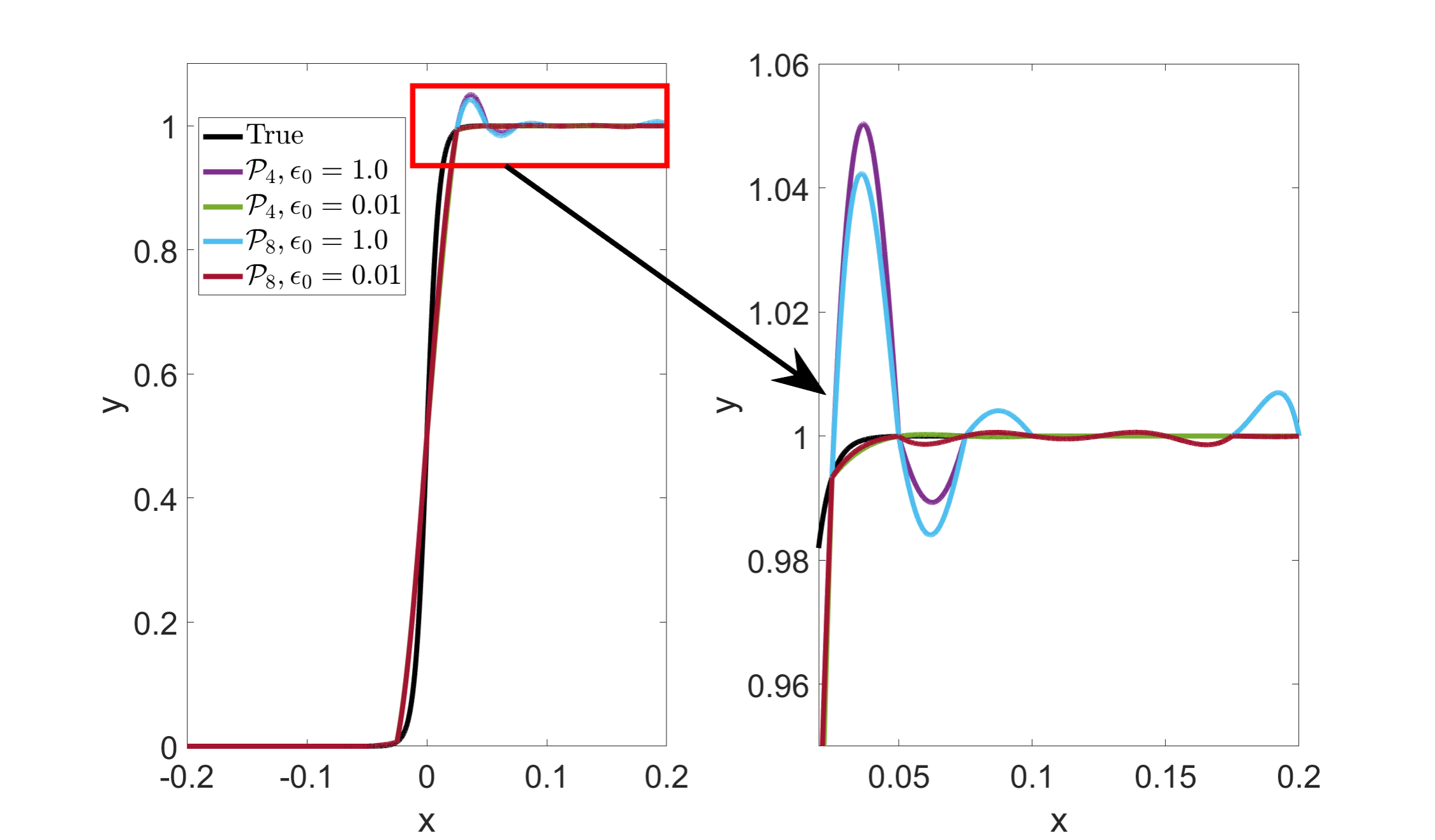}
            \subcaption{PPI with $d=4$ and $d=8$.}
            \label{subfig:heaviside2}
        \end{subfigure}
        \caption{Approximation of $f_2(x)$ from $N=17$ uniformly spaced points with different interpolation methods. 
        The top row (Figure \ref{subfig:heaviside1}) shows approximation results using PCHIP, MQSI, and DBI.
        The bottom row (Figure \ref{subfig:heaviside2}) shows approximation results using PPI with $d=4, 8$ and $\epsilon_{0}=1.0, 0.01$.
        An enlarged version of the region in the red rectangle is shown on the right of each row.
        The value of $\epsilon_{1}$ is set to $1.0$.}
        \label{fig:heaviside}
    \end{figure}

      Figure \ref{fig:heaviside} shows approximation plots of $f_{2}(x)$ using $N=17$ uniformly spaced points with different values of $\epsilon_{0}$ and $\epsilon_{1}=1$.
      The target polynomial degree is set to $d=8$.
      For $\epsilon_{0} =1$, we observe oscillations, as shown in the right part of Figure \ref{fig:heaviside}.
      As $\epsilon_{0}$ decreases, the oscillations decrease. 
      For $\epsilon_{0} \leq 0.01$, the errors and oscillations are negligible compared to errors in the region with the steep gradient. 
      The oscillations are completely removed for $\epsilon_{0} =0.0$. 
      \begin{table}[H]
        \centering
        \begin{tabular}{ c c c c c c c c c c c c c}
          \hline
          \hline
          $N$   && PCHIP  && MQSI  && \multicolumn{3}{c}{DBI} && \multicolumn{3}{c}{PPI}   \\

                &&$\mathcal{P}_{3}$ && $\mathcal{P}_{5}$ && $\mathcal{P}_{3}$ & $\mathcal{P}_{4}$ & $\mathcal{P}_{8}$
                                    &&$\mathcal{P}_{3}$ & $\mathcal{P}_{4}$ & $\mathcal{P}_{8}$ \\
          \hline
          17       && 2.02E-2  && 1.92E-2 &&  2.41E-2  &  2.41E-2  &  2.08E-2  &&  2.41E-2  &  2.41E-2  &  2.08E-2   \\
          33       && 3.38E-3  && 3.72E-3 &&  4.89E-3  &  4.86E-3  &  3.59E-3  &&  4.90E-3  &  4.86E-3  &  3.57E-3   \\
          65       && 3.59E-4  && 1.61E-3 &&  4.17E-4  &  1.89E-4  &  1.47E-4  &&  4.17E-4  &  1.89E-4  &  1.47E-4   \\
          129      && 4.21E-5  && 1.71E-4 &&  3.09E-5  &  1.55E-5  &  1.70E-6  &&  3.09E-5  &  1.55E-5  &  1.70E-6   \\
          257      && 5.12E-6  && 1.75E-5 &&  2.04E-6  &  5.31E-7  &  5.22E-9  &&  2.04E-6  &  5.31E-7  &  5.22E-9   \\
          \hline
          \hline
        \end{tabular}                                                                                  
        \caption{$L^2$-errors when using the PCHIP, MQSI , DBI, and PPI  methods to approximate the function $f_{2}(x)$.
                 $N$ represents the number of input points used to build the approximation.
                 The parameters $\epsilon_{0}$, $\epsilon_{1}$, and $st$ are set to $0.01, 1,$ and $3$, respectively.}
        \label{tab:f2}
      \end{table}
      
    \subsection{ Example III: 1D Discontinuous Function }
      This example uses a modified version of a function introduced by Tadmor and Tanner \cite{Tadmor2002} and used by Berzins \cite{Berzins} in the context of bounded interpolation methods.
      The value one is added to the original function in ~\cite{Tadmor2002} to ensure that the function is positive over the interval [-1,1]. 
      The modified function is defined as
      \begin{equation}\label{eq:GT}
        f_{3}(x) =
        \begin{cases}
          1+\frac{2 e^{2\pi x} -1 -e^{\pi}}{e^{\pi}-1}, \quad x \in [-1, -0.5)\\\\
          1- sin \big(\frac{2\pi x}{3}+ \frac{\pi}{3}\big), \quad x \in [-0.5, 1].
        \end{cases}
      \end{equation}

      Table \ref{tab:f3} shows $L^{2}$-error norms when using the PCHIP, MQSI, DBI, and PPI methods to approximate the function in Equation (\ref{eq:GT}).
      Approximating $f_{3}(x)$ is challenging because $f_{3}(x)$ is a piecewise function with a discontinuity at $x=0.5$. 
      The global error is dominated by the local errors around the discontinuity.
      The PCHIP, MQSI, DBI, and PPI approximation results are comparable. 
      Increasing the target polynomial degree does not decrease the $L^{2}$-error norms.
      The approximations in the smooth regions improve as we increase the target polynomial degree, but the global error is dominated by the error around the discontinuity.
      The error around the discontinuity does not decrease with higher polynomial degrees.
      
\begin{figure}
        \centering
        \begin{subfigure}{0.9 \textwidth}
            \centering
            \includegraphics[scale=0.35]{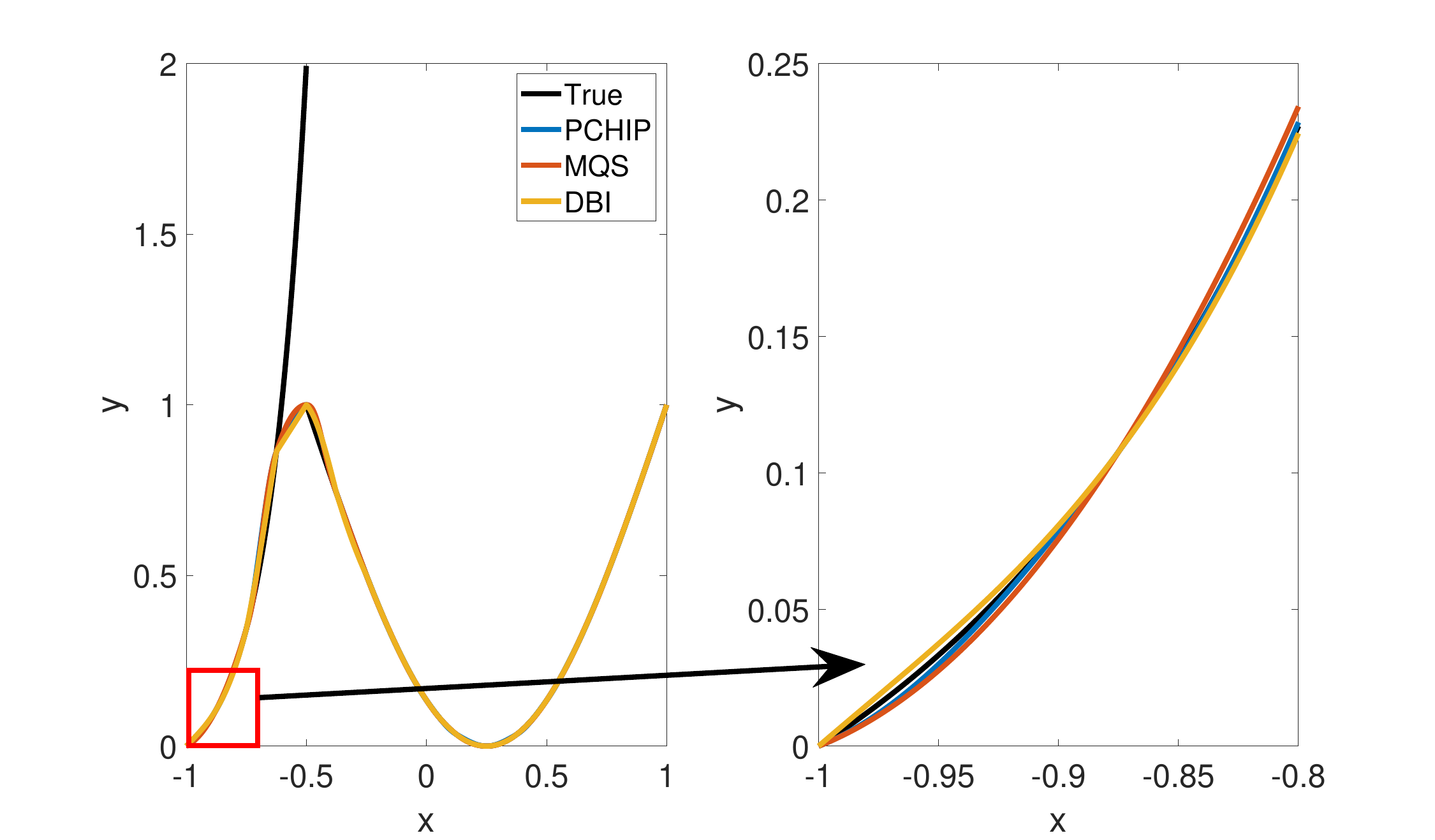}
            \subcaption{PCHIP, MQSI, and DBI with $d=8$.}
            \label{subfig:GelbT1}
        \end{subfigure}
        \begin{subfigure}{0.9 \textwidth}
            \centering
            \includegraphics[scale=0.35]{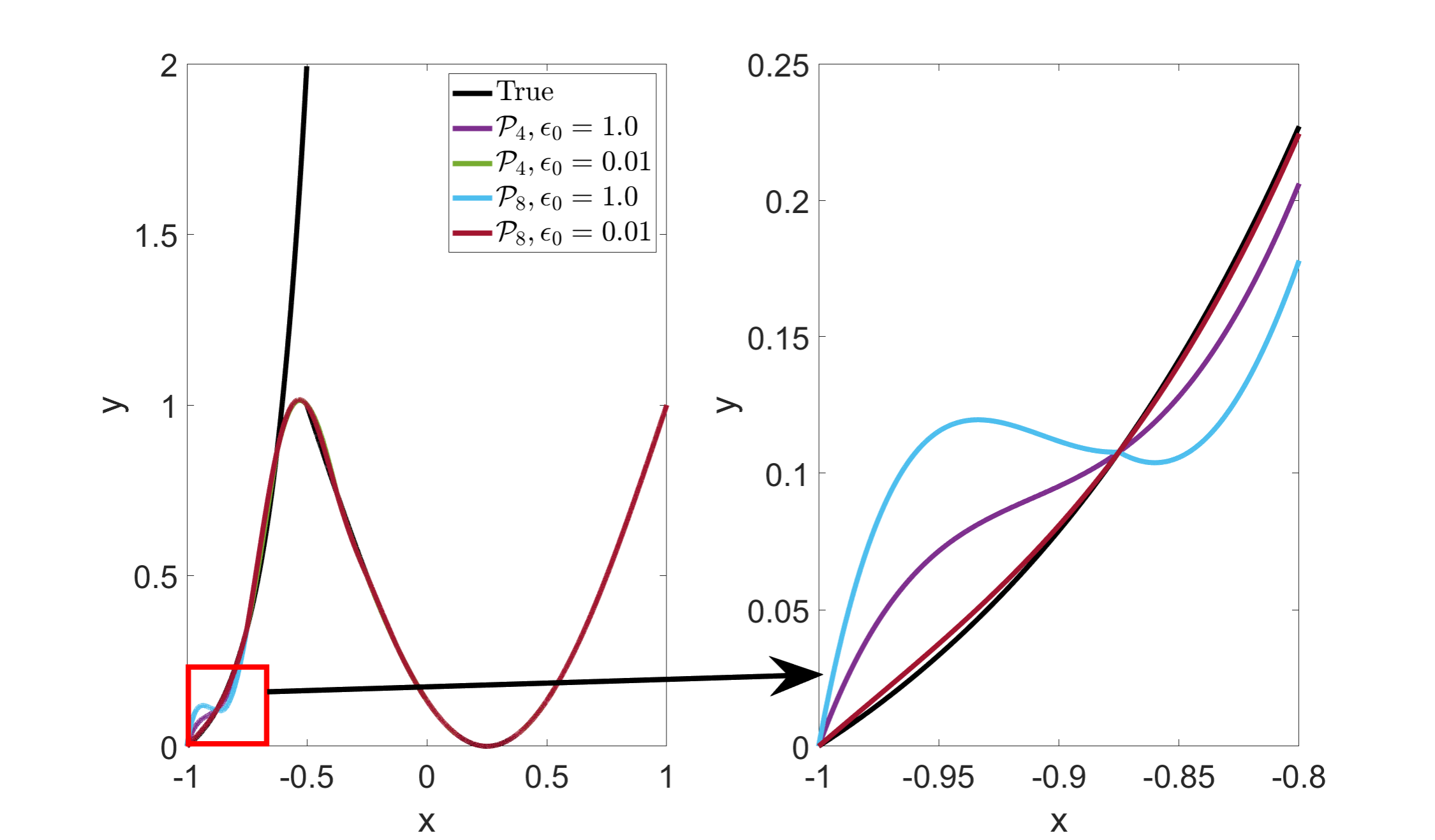}
            \subcaption{PPI with $d=4$ and $d=8$.}
            \label{subfig:GelbT2}
        \end{subfigure}
        \caption{Approximation of $f_3(x)$ from $N=17$ uniformly spaced points with different interpolation methods. 
        The top row (Figure \ref{subfig:GelbT1}) shows approximation results using PCHIP, MQSI, and DBI.
        The bottom row (Figure \ref{subfig:GelbT2}) shows approximation results using PPI with $d=4, 8$ and $\epsilon_{0}=1.0, 0.01$.
        An enlarged version of the region in the red rectangle is shown on the right of each row.
        The value of $\epsilon_{1}$ is set to $1.0$.}
        \label{fig:GelbT}
    \end{figure}

      Figure \ref{fig:GelbT} shows approximation plots of $f_{3}(x)$ using $N=17$ uniformly spaced points with different values of $\epsilon_{0}$.
      The target polynomial degree is set to $d=4, 8$.
      The bottom right part of Figure \ref{fig:GelbT} shows oscillations at the left boundary for $\epsilon_{0} =1$.
      The oscillations are removed for $\epsilon_{0} \leq 0.1$.
      As expected, all the interpolation methods have difficulties approximating the function around the discontinuity, as shown in Figure \ref{fig:GelbT}. 
      \begin{table}[H]
        \centering
        \begin{tabular}{ c c c c c c c c c c c c c c}
          \hline
          \hline
          $N$   && PCHIP  &&  MQSI  && \multicolumn{3}{c}{DBI} && \multicolumn{3}{c}{PPI}   \\

                &&$\mathcal{P}_{3}$ &&  $\mathcal{P}_{5}$  && $\mathcal{P}_{3}$ & $\mathcal{P}_{4}$ & $\mathcal{P}_{8}$
                                    &&$\mathcal{P}_{3}$ & $\mathcal{P}_{4}$ & $\mathcal{P}_{8}$ \\
          \hline
          17       && 1.77E-1  && 1.76E-1 &&  1.82E-1  &  1.83E-1  &  1.82E-1  &&  1.73E-1  &  1.72E-1  &  1.70E-1   \\
          33       && 1.39E-1  && 1.41E-1 &&  1.35E-1  &  1.39E-1  &  1.36E-1  &&  1.35E-1  &  1.39E-1  &  1.36E-1   \\
          65       && 1.03E-1  && 1.06E-1 &&  9.95E-2  &  1.04E-1  &  1.02E-1  &&  9.95E-2  &  1.04E-1  &  1.02E-1   \\
          129      && 7.42E-2  && 7.63E-2 &&  7.12E-2  &  7.54E-2  &  7.35E-2  &&  7.15E-2  &  7.55E-2  &  7.38E-2   \\
          257      && 5.28E-2  && 5.43E-2 &&  5.06E-2  &  5.38E-2  &  5.24E-2  &&  5.07E-2  &  5.39E-2  &  5.26E-2   \\
          \hline
          \hline
        \end{tabular}                                                                                  
        \caption{$L^2$-errors when using the PCHIP, MQSI DBI, and PPI  methods to approximate the function $f_{3}(x)$.
                 $N$ represents the number of input points used to build the approximation.
                 The parameters $\epsilon_{0}$, $\epsilon_{1}$, and $st$ are set to $0.01, 1,$ and $3$, respectively.}
        \label{tab:f3}
      \end{table}

%\newpage

    \subsection{Example IV: 2D Modified Runge Function}
      This example extends the previously modified 1D Runge function to 2D as follows:
      \begin{equation}\label{eq:f4}
        f_{4}(x, y) = \frac{0.1}{0.1+25(x^{2}+y^{2})}, \quad x, y \in [-1,1].
      \end{equation}

      Table \ref{tab:f4} shows $L^{2}$-error norms when using the PCHIP, MQSI, DBI, and PPI methods to approximate the 2D modified Runge function in Equation (\ref{eq:f4}).
      In this example, as in Example I, the MQSI and PCHIP methods have comparable errors.
      The approaches used to approximate the derivatives and construct the quintic splines are not more accurate than the approximation obtained from PCHIP.
      The DBI and PPI methods with $d=3$ lead to better approximation results compared to the PCHIP and MQSI. 
      As the target polynomial degree $d$ increases, the approximation errors from PPI decrease much faster than from DBI. 
      The relaxed nature of the PPI methods allows for higher degree polynomials compared to DBI.
      The bounds for data-boundedness, which are based on Equation (\ref{eq:uminumax}) with $\Delta_{min} = \Delta_{max} =0.0$, are more restrictive than positivity for which $\Delta_{min} > 0.0$ and $ \Delta_{max} >0.0$.
      In addition, the approximation does not lead to oscillations for $\epsilon_{0}$ and $\epsilon_{1} \in [0,1]$

      \begin{table}[H]
        \centering
        \begin{tabular}{ c c c c c c c c c c c c c}
          \hline
          \hline
          $N^{2}$      && PCHIP   &&  MQSI  && \multicolumn{3}{c}{DBI} && \multicolumn{3}{c}{PPI}   \\

                   &&$\mathcal{P}_{3}$ && $\mathcal{P}_{5}$ &&$\mathcal{P}_{3}$ & $\mathcal{P}_{4}$ & $\mathcal{P}_{8}$
                                    &&$\mathcal{P}_{3}$ & $\mathcal{P}_{4}$ & $\mathcal{P}_{8}$ \\
          \hline
          17       && 1.76E-2  && 1.67E-2  &&  2.12E-2  &  9.09E-3  &  1.91E-2  &&  2.12E-2  &  9.09E-3  &  1.91E-2   \\
          33       && 2.05E-3  && 2.84E-3  &&  2.45E-3  &  4.61E-3  &  1.25E-3  &&  2.45E-3  &  4.61E-3  &  1.24E-3   \\
          65       && 1.05E-3  && 8.01E-4  &&  8.59E-4  &  9.33E-4  &  4.99E-4  &&  8.59E-4  &  9.33E-4  &  3.51E-4   \\
          129      && 2.23E-4  && 2.57E-4  &&  7.47E-5  &  4.76E-5  &  4.12E-5  &&  7.47E-5  &  4.64E-5  &  7.16E-6   \\
          257      && 4.19E-5  && 3.53E-5  &&  5.05E-6  &  4.20E-6  &  3.80E-6  &&  5.05E-6  &  1.62E-6  &  2.91E-8   \\
          \hline
          \hline
        \end{tabular}                                                                                  
        \caption{$L^2$-errors when using the PCHIP, MQSI, DBI, and PPI  methods to approximate the function $f_{4}(x,y)$.
                 $N$ represents the number of input points used to build the approximation.
                 The parameters $\epsilon_{0}$, $\epsilon_{1}$, and $st$ are set to $0.01, 1,$ and $3$, respectively.} 
        \label{tab:f4}
      \end{table}

    \subsection{Example V: 2D Logistic Function}
      The test case extends the 1D logistic (smoothed Heaviside) function from Example II to 2D.
      \begin{equation}\label{eq:Heaviside2}
        f_{5}(x,y) = \frac{1}{1+e^{-\sqrt{2}k(x+y)}}, \quad x, y \in [-0.2, 0.2].
      \end{equation}
      The function $f_{5}(x,y)$ is challenging because of the large gradient at $y=-x$.
      Approximating $f_{5}(x,y)$ with a standard polynomial interpolation leads to oscillations and negative values that violate the desired property of positivity.

\begin{figure}
        \centering
        \begin{subfigure}{0.48 \textwidth}
            \centering
            \includegraphics[scale=0.48]{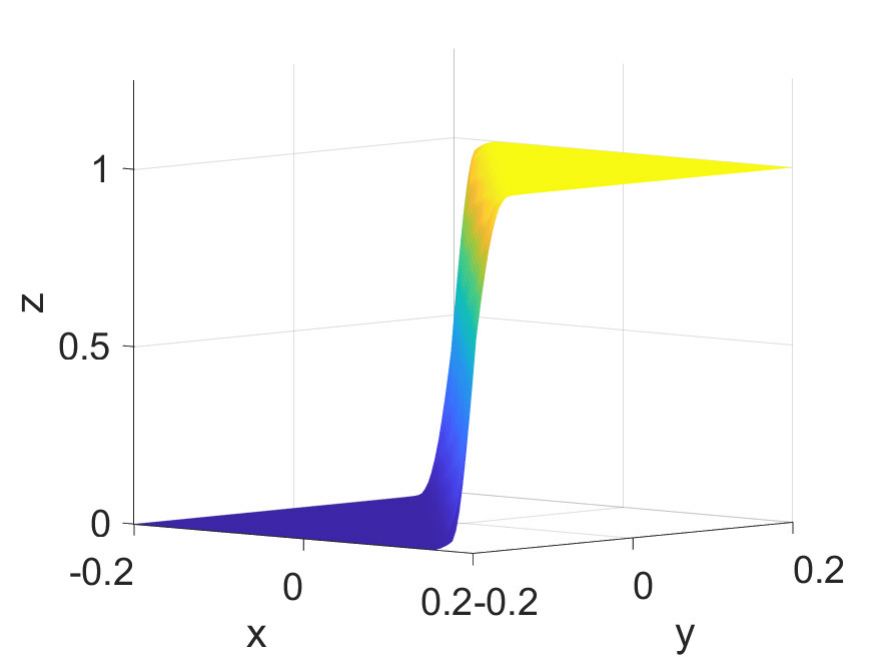}
            \subcaption{PCHIP}
            \label{subfig:heaviside2DPCHIP}
        \end{subfigure}
        \begin{subfigure}{0.48 \textwidth}
            \centering
            \includegraphics[scale=0.48]{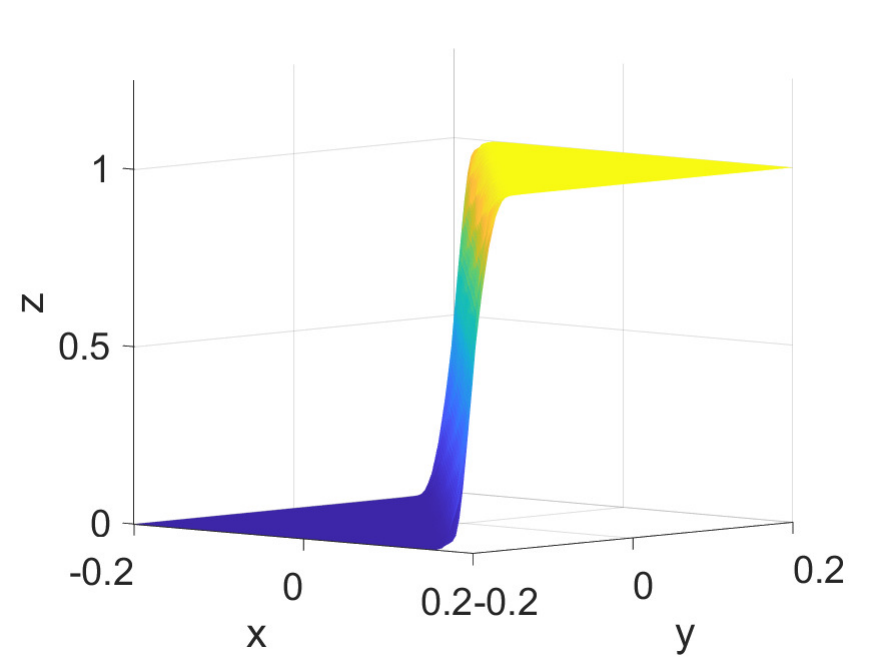}
            \subcaption{MQSI}
            \label{subfig:heaviside2DMQSI}
        \end{subfigure}
        \begin{subfigure}{0.48 \textwidth}
            \centering
            \includegraphics[scale=0.48]{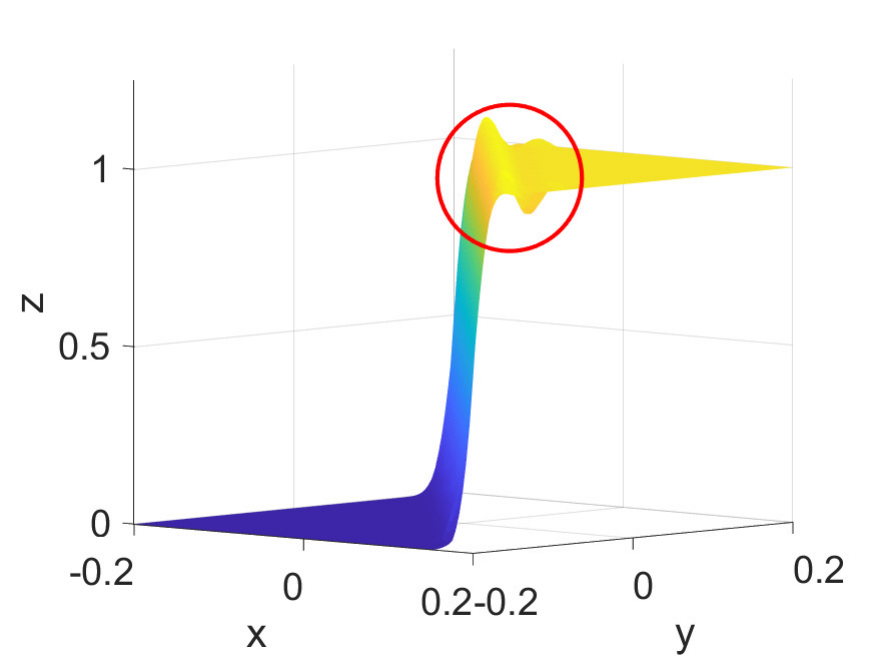}
            \subcaption{PPI with $\mathcal{P}_{4}$, $\epsilon_{0}=\epsilon_{1}=1.0$}
            \label{subfig:heaviside2DP4_1}
        \end{subfigure}
        \begin{subfigure}{0.48 \textwidth}
            \centering
            \includegraphics[scale=0.48]{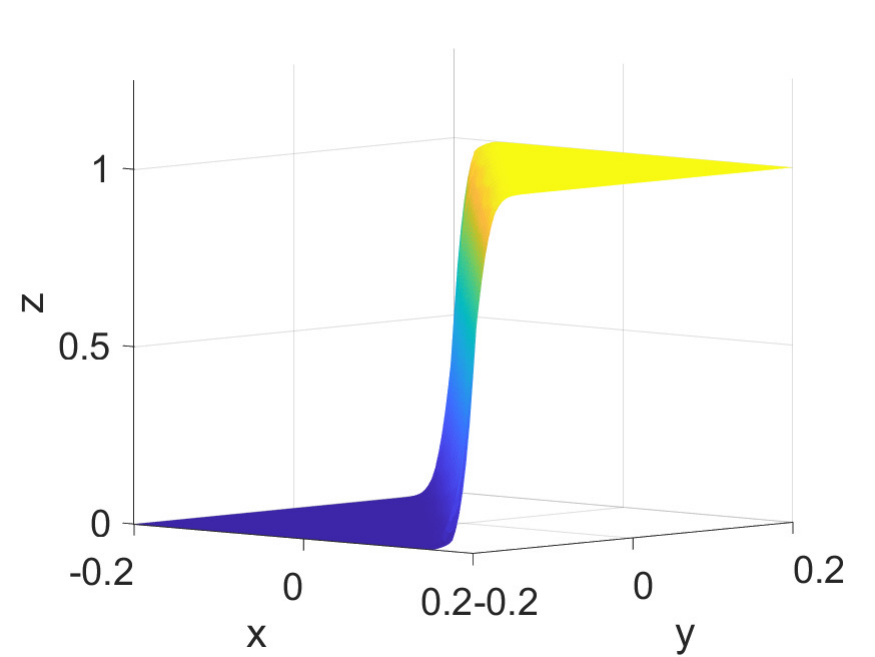}
            \subcaption{PPI with $\mathcal{P}_{4}$, $\epsilon_{0}=\epsilon_{1}=10^{-4}$}
            \label{subfig:heaviside2DP4_2}
        \end{subfigure}
        \begin{subfigure}{0.48 \textwidth}
            \centering
            \includegraphics[scale=0.48]{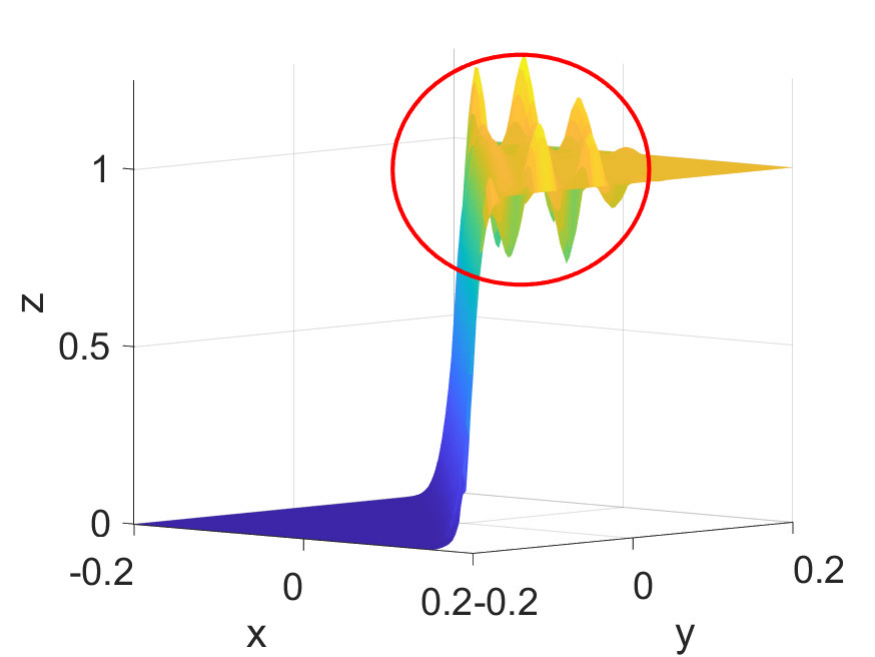}
            \subcaption{PPI with $\mathcal{P}_{8}$, $\epsilon_{0}=\epsilon_{1}=1.0$}
            \label{subfig:heaviside2DP8_1}
        \end{subfigure}
        \begin{subfigure}{0.48 \textwidth}
            \centering
            \includegraphics[scale=0.48]{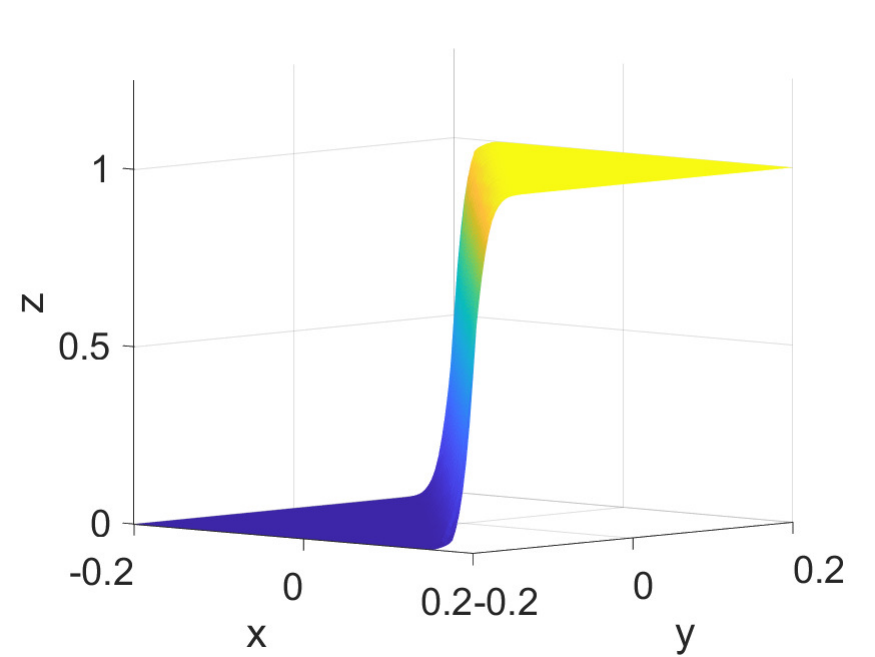}
            \subcaption{PPI with $\mathcal{P}_{8}$, $\epsilon_{0}=\epsilon_{1}=10^{-4}$}
            \label{subfig:heaviside2DP8_2}
        \end{subfigure}
        \caption{Approximation of $f_5(x,y)$ from $N \times N =17^2$ uniformly spaced points with different interpolation methods. 
        The parameter $st$ is set to $2$.
        The red ellipses highlight examples regions with large oscillations.
        The surfaces in the top row Figures \ref{subfig:heaviside2DPCHIP} and \ref{subfig:heaviside2DMQSI} are obtained using PCHIP and MQSI.
        Figures \ref{subfig:heaviside2DP4_1} ($\mathcal{P}_{4}$, $\epsilon_{0}=\epsilon_{1}=1.0$), \ref{subfig:heaviside2DP4_2} ($\mathcal{P}_{4}$, $\epsilon_{0}=\epsilon_{1}=10^{-4}$), \ref{subfig:heaviside2DP8_1} ($\mathcal{P}_{8}$, $\epsilon_{0}=\epsilon_{1}=1.0$), and \ref{subfig:heaviside2DP8_2} ($\mathcal{P}_{8}$, $\epsilon_{0}=\epsilon_{1}=10^{-4}$) are obtained using PPI. }
        \label{fig:heaviside2D}
    \end{figure}

      Table \ref{tab:f5} shows $L^{2}$-error norms when using the PCHIP, MQSI, DBI, and PPI methods to approximate the 2D logistic function $f_{5}(x, y)$ defined in Equation (\ref{eq:Heaviside2}). 
      The errors from MQSI are larger compared to the other methods.
      In this example, the approximated derivatives and quintic splines from MQSI are less accurate than the approximations from PCHIP. 
      The DBI and PPI methods lead to better approximation results compared to the PCHIP and MQSI approaches.
      Increasing the target polynomial degree improves the approximations for DBI and PPI, as shown in Table \ref{tab:f5}.
      The global error is dominated by the local around the steep gradients at $y=-x$.
      For $N^{2}=33^{2}$ points and above with $\epsilon_{0}=0.01$, $\epsilon_{1}=1.0$, and $st=3$, the DBI and PPI choose the identical stencils for the intervals around the steep gradient. 
      Therefore, the global error, which is dominated by the local error around the steep gradient, is the same for both DBI and PPI.
      However, for $N^{2}=17^{2}$ points with $\epsilon_{0}=0.01$, $\epsilon_{1}=1.0$, and $st=3$, the DBI and PPI choose different stencils that lead to different errors, as indicated in the row with $N^{2}=17^{2}$ of Table \ref{tab:f5}.
      Figure \ref{fig:heaviside2D} shows approximation plots of $f_{5}(x,y)$ using $N\times N=17\times 17$ uniformly spaced points with PCHIP, MQSI, and PPI.
      The approximations in Figures \ref{subfig:heaviside2DP4_1} and \ref{subfig:heaviside2DP8_1} are obtained using PPI with $\epsilon_{0}=\epsilon_{1}=1$.
      The PPI method with $\epsilon_{0}=\epsilon_{1}=10^{-4}$ is used for the approximations in Figures \ref{subfig:heaviside2DP4_2} and \ref{subfig:heaviside2DP8_2}.
      The target polynomial degrees for the second and third rows are set to $d=4$ and $d=8$, respectively.
      The oscillations increase when going from $\mathcal{P}_{4}$ to $\mathcal{P}_{8}$ with $\epsilon_{1}=\epsilon_{2} =1.0$.
      For $\epsilon_{0} =\epsilon_{1}=10^{-4}$, the oscillations are significantly reduced, and the approximation is closer to the target solution.

      \begin{table}[H]
        \centering
        \begin{tabular}{ c c c c c c c c c c c c c}
          \hline
          \hline
          $N^{2}$  && PCHIP &&  MQSI   && \multicolumn{3}{c}{DBI} && \multicolumn{3}{c}{PPI}   \\
                   &&$\mathcal{P}_{3}$ && $\mathcal{P}_{5}$ && $\mathcal{P}_{3}$ & $\mathcal{P}_{4}$ & $\mathcal{P}_{8}$
                                    &&$\mathcal{P}_{3}$ & $\mathcal{P}_{4}$ & $\mathcal{P}_{8}$ \\
          \hline
          17       && 8.07E-3  && 7.12E-3  &&  1.05E-2  &  9.79E-3  &  8.18E-3  &&  1.05E-2  &  9.77E-3  &  8.61E-3   \\
          33       && 1.26E-3  && 2.62E-3  &&  1.67E-3  &  1.36E-3  &  1.06E-3  &&  1.64E-3  &  1.30E-3  &  8.87E-4   \\
          65       && 1.44E-4  && 6.35E-4  &&  1.58E-4  &  8.84E-5  &  4.89E-5  &&  1.58E-4  &  8.84E-5  &  5.01E-5   \\
          129      && 1.63E-5  && 6.51E-5  &&  1.13E-5  &  3.07E-6  &  2.64E-7  &&  1.13E-5  &  3.07E-6  &  2.64E-7   \\
          257      && 1.94E-6  && 6.74E-6  &&  7.29E-7  &  1.02E-7  &  5.39E-10  &&  7.29E-7  &  1.02E-7  &  5.39E-10   \\
          \hline
          \hline
        \end{tabular}                                                                                  
        \caption{$L^2$-errors when using the PCHIP, MQSI, DBI, and PPI  methods to approximate the function $f_{5}(x,y)$.
                 $N$ represents the number of input points used to build the approximation.
                 The parameters $\epsilon_{0}$, $\epsilon_{1}$, and $st$ are set to $0.01, 1,$ and $3$, respectively.}
        \label{tab:f5}
      \end{table}

    \subsection{Example VI: \texorpdfstring{$C^{0}$}{TEXT}-continuous Surface Function}
      This example is based on a 2D function used to study positive and monotonic splines \cite{CHAN2001135, 10.1007/11537908_20}.
      The function is defined as follows:
      \begin{equation}\label{eq:T2}
        f_{6}(x, y) =
        \begin{cases}
          2(y-x) & \textrm{ if } 0 \leq y-x \leq 0.5 \\
          1      & \textrm{ if } y-x \geq 0.5 \\
          cos \bigg( 4\pi \sqrt{(x-1.5)^{2}+(y-0.5)^{2}}\bigg) & \textrm{ if }(x-1.5)^{2}+(y-0.5)^{2} \leq \frac{1}{16} \\
          0      & otherwise.
        \end{cases}
      \end{equation}
      The function $f_{6}(x,y)$ is challenging because it is only $\mathbb{C}^{0}$-continuous at various locations. 

      Table \ref{tab:f6} shows $L^{2}$-error norms when using the PCHIP, MQSI, DBI, and PPI methods to approximate the 2D function $f_{6}(x, y)$ defined in Equation (\ref{eq:T2}).
      The PCHIP, MQSI, DBI, and PPI methods lead to comparable $L^{2}-$error norms.
      Increasing the target polynomial degree does not significantly improve the approximation for DBI and PPI, as shown in Table \ref{tab:f6}.
      The global error is dominated by the local around the $\mathbb{C}^{0}$.
      The approximation for both DBI and PPI can be improved by using an underlying mesh that better captures the $\mathbb{C}^{0}$-continuity.

\begin{figure}
        \centering
        \begin{subfigure}{0.48 \textwidth}
            \centering
            \includegraphics[scale=0.48]{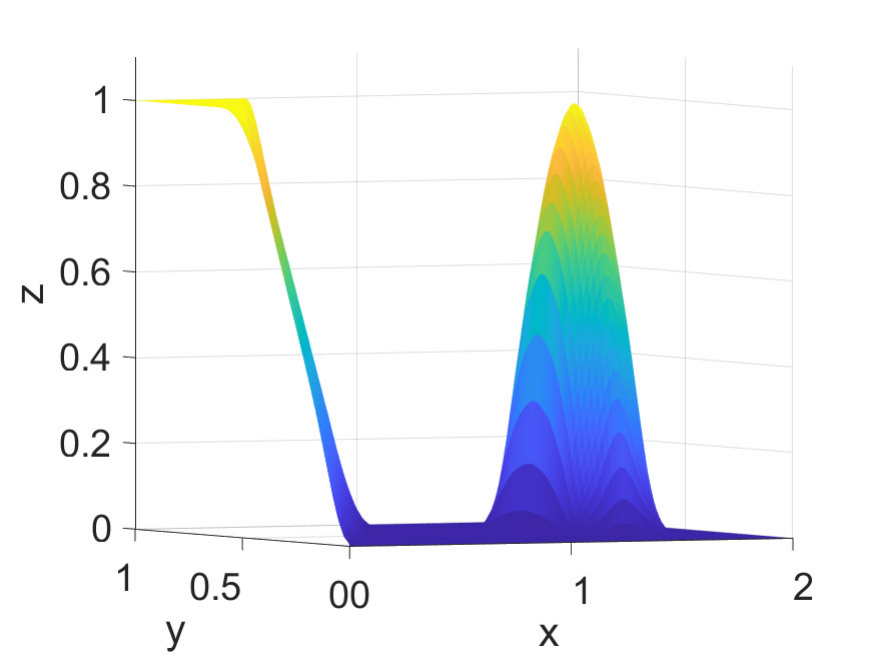}
            \subcaption{PCHIP}
            \label{subfig:surface2PCHIP}
        \end{subfigure}
        \begin{subfigure}{0.48 \textwidth}
            \centering
            \includegraphics[scale=0.48]{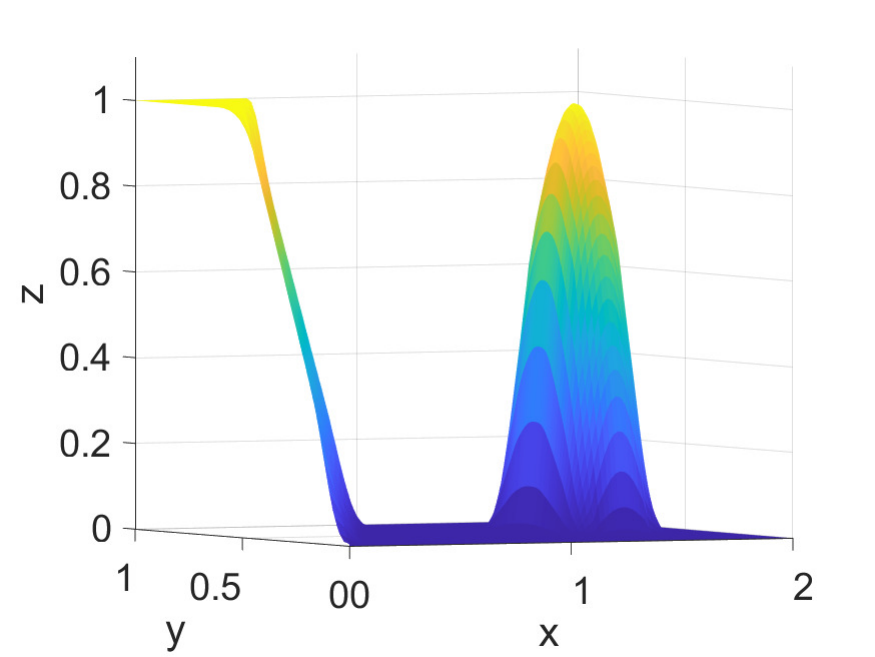}
            \subcaption{MQSI}
            \label{subfig:surface2MQSI}
        \end{subfigure}
        \begin{subfigure}{0.48 \textwidth}
            \centering
            \includegraphics[scale=0.48]{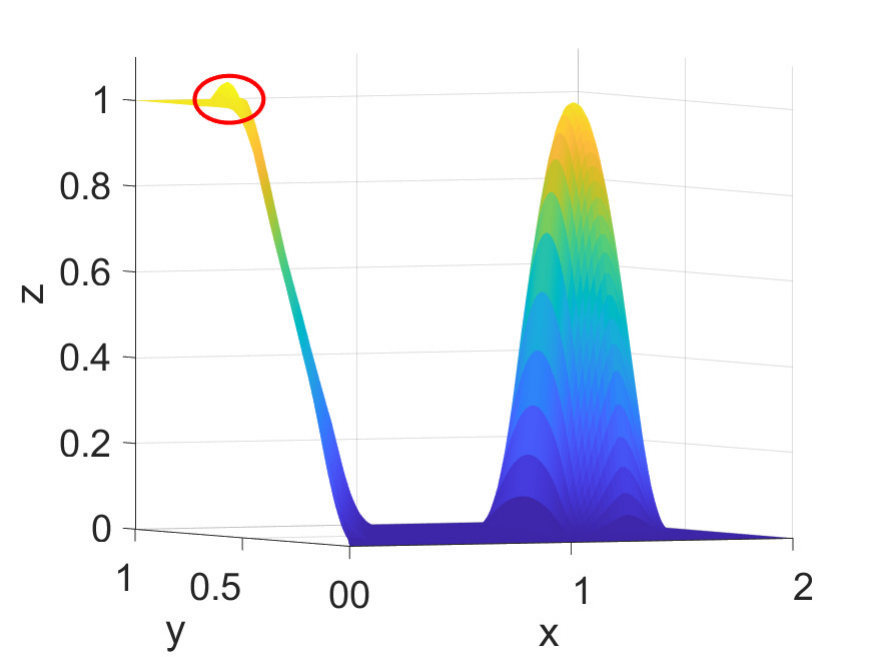}
            \subcaption{PPI with $\mathcal{P}_{4}$, $\epsilon_{0}=\epsilon_{1}=1.0$}
            \label{subfig:surface2P4_1}
        \end{subfigure}
        \begin{subfigure}{0.48 \textwidth}
            \centering
            \includegraphics[scale=0.48]{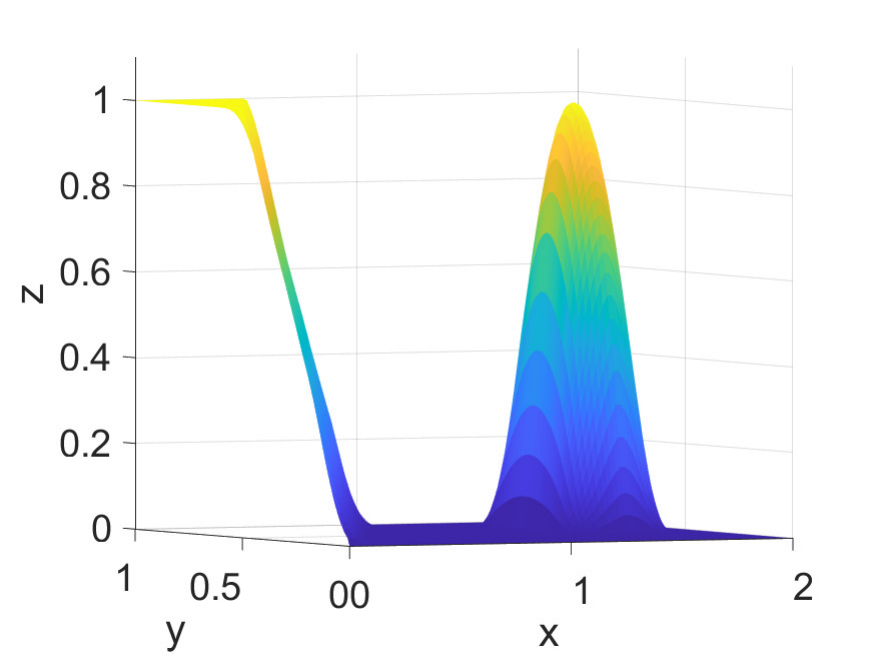}
            \subcaption{PPI with $\mathcal{P}_{4}$, $\epsilon_{0}=\epsilon_{1}=10^{-4}$}
            \label{subfig:surface2P4_2}
        \end{subfigure}
        \begin{subfigure}{0.48 \textwidth}
            \centering
            \includegraphics[scale=0.48]{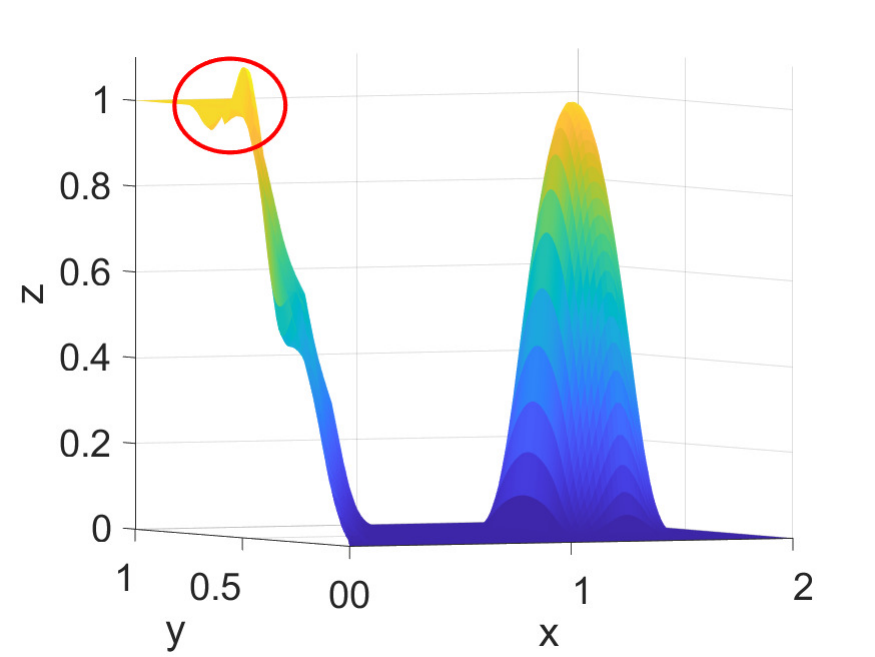}
            \subcaption{PPI with $\mathcal{P}_{8}$, $\epsilon_{0}=\epsilon_{1}=1.0$}
            \label{subfig:surface2P8_1}
        \end{subfigure}
        \begin{subfigure}{0.48 \textwidth}
            \centering
            \includegraphics[scale=0.48]{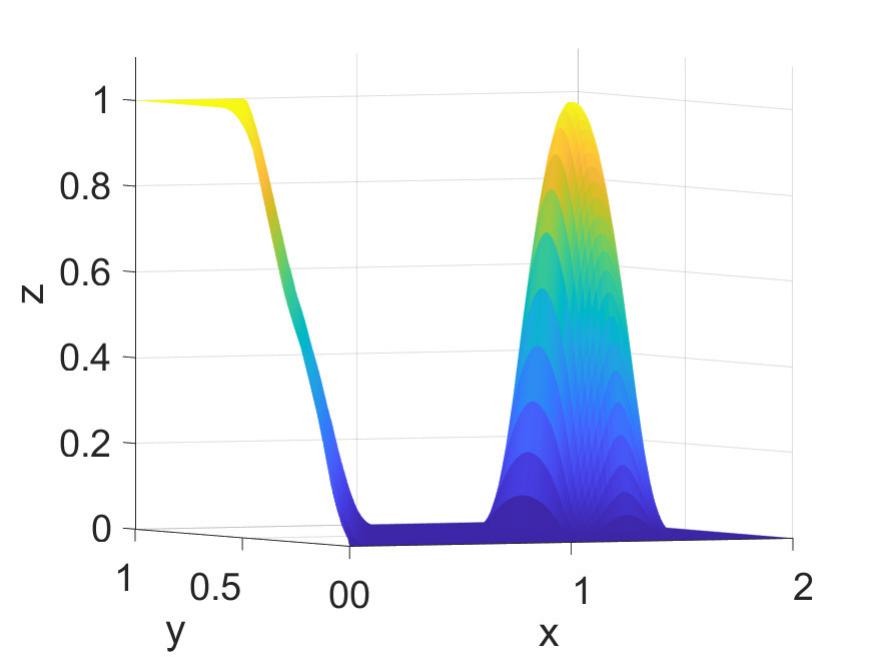}
            \subcaption{PPI with $\mathcal{P}_{8}$, $\epsilon_{0}=\epsilon_{1}=10^{-4}$}
            \label{subfig:surface2P8_2}
        \end{subfigure}
        \caption{Approximation of $f_5(x,y)$ from $N \times N =17^2$ uniformly spaced points with different interpolation methods. 
        The parameter $st$ is set to $2$.
        The red ellipses highlight examples regions with large oscillations.
        The surfaces in the top row Figures \ref{subfig:surface2PCHIP} and \ref{subfig:surface2MQSI} are obtained using PCHIP and MQSI.
        Figures \ref{subfig:surface2P4_1} ($\mathcal{P}_{4}$, $\epsilon_{0}=\epsilon_{1}=1.0$), \ref{subfig:surface2P4_2} ($\mathcal{P}_{4}$, $\epsilon_{0}=\epsilon_{1}=10^{-4}$), \ref{subfig:surface2P8_1} ($\mathcal{P}_{8}$, $\epsilon_{0}=\epsilon_{1}=1.0$), and \ref{subfig:surface2P8_2} ($\mathcal{P}_{8}$, $\epsilon_{0}=\epsilon_{1}=10^{-4}$) are obtained using PPI.}
        \label{fig:surface2}
    \end{figure}

      Figure \ref{fig:surface2} shows approximation plots of $f_{6}(x,y)$ using $N\times N=17\times 17$ uniformly spaced points with PCHIP, MQSI, and PPI.
      The left and right approximation plots show approximated solutions using the PPI method.
      For the left plot using PPI, $\epsilon_{0}=\epsilon_{1}=1.0$, and for the right plot using PPI, $\epsilon_{0}=\epsilon_{1}=10^{-4}$
      The target polynomial degree is set to $d=4$ and $d=8$ for the second and third rows, respectively.
      For $\epsilon_{0}=\epsilon_{1}=1.0$, the oscillations increase as the target polynomial degree increases from $d=4$ to $d=8$. 
      The oscillations observed for $\epsilon_{1}=\epsilon_{0}=1$ are removed for small values of $\epsilon_{0}$ and $\epsilon_{1}$, shown in Figures \ref{subfig:surface2P4_2} and \ref{subfig:surface2P8_2}. 

      \begin{table}[H]
        \centering
        \begin{tabular}{ c c c c c c c c c c c c c}
          \hline
          \hline
          $N^{2}$  && PCHIP && MQSI  && \multicolumn{3}{c}{DBI} && \multicolumn{3}{c}{PPI }   \\
                   &&$\mathcal{P}_{3}$ &&  $\mathcal{P}_{5}$  && $\mathcal{P}_{3}$ & $\mathcal{P}_{4}$ & $\mathcal{P}_{8}$
                                    &&$\mathcal{P}_{3}$ & $\mathcal{P}_{4}$ & $\mathcal{P}_{8}$ \\
          \hline
          17       && 1.91E-2  && 2.19E-2  &&  1.72E-2  &  1.69E-2  &  1.63E-2  &&  1.72E-2  &  1.68E-2  &  1.59E-2   \\
          33       && 6.92E-3  && 7.44E-3  &&  6.16E-3  &  5.81E-3  &  5.88E-3  &&  6.16E-3  &  5.80E-3  &  5.87E-3   \\
          65       && 2.47E-3  && 2.45E-3  &&  2.24E-3  &  2.14E-3  &  2.11E-3  &&  2.24E-3  &  2.14E-3  &  2.11E-3   \\
          129      && 8.99E-4  && 8.69E-4  &&  8.21E-4  &  7.77E-4  &  7.63E-4  &&  8.20E-4  &  7.77E-4  &  7.63E-4   \\
          257      && 3.23E-4  && 3.04E-4  &&  2.97E-4  &  2.81E-4  &  2.76E-4  &&  2.96E-4  &  2.81E-4  &  2.76E-4   \\
          \hline                                                                         
          \hline                                                              
        \end{tabular}                                                                                  
        \caption{$L^2$-errors when using the PCHIP, MQSI, DBI, and PPI  methods to approximate the function $f_{6}(x,y)$.
                 $N$ represents the number of input points used to build the approximation.
                 The parameters $\epsilon_{0}$, $\epsilon_{1}$, and $st$ are set to $0.01, 1,$ and $3$, respectively.}
        \label{tab:f6}
      \end{table}
   
  \section{Mapping Error Analysis and Examples}
\label{sec:error-mapping-examples}
    This section provides an analysis of the mapping error between two different meshes in the context of time-dependent PDEs when DBI and  PPI are employed to interpolate data values between the meshes.
    In addition to the error analysis, the Runge, tropical warm pool international cloud experiment (TWP-ICE), and Barbados oceanographic and meteorological experiment (BOMEX) examples are used to evaluate the use of the PPI and DBI to map data values between two meshes.
    The Runge and TWP-ICE examples use meshes that emulate the advection and reaction meshes used in NEPTUNE.
    These meshes are constructed by linearly scaling the NEPTUNE vertical mesh points to the desired interval for the Runge and TWP-ICE examples. 
    In the BOMEX example, the advection mesh is composed of uniformly spaced points, and the reaction mesh is constructed using the mid-point of each interval from the advection mesh.
    %$In Examples \ref{subsec:runge} and \ref{subsec:twp-ice}, we consider the max error instead of $L^{2}$-norm error because the global error is dominated by the local errors in a few locations around the large gradients.  
   
  \subsection{Mapping Error Analysis}
  \label{subsec:error}
  In addition to the development and study of the DBI and PPI methods, it is important to provide some insight into the behavior of the mapping error in the context of time-dependent PDEs.
  An example of a time-dependent problem where a positivity-preserving mapping is required is the US Navy Environmental Prediction System Using the NUMA Core (NEPTUNE) \cite{doyle:NEPTUNE}. 
  NEPTUNE is a next-generation global NWP system being developed at the Naval Research Laboratory (NRL) and the Naval Postgraduate School (NPS). 
  In NEPTUNE, the physics and dynamics are calculated using different meshes and require mapping the solution values between both meshes. 
  NEPTUNE uses nonuniform structured meshes that have vertical columns with nonuniformly spaced points inside each column. 
  The mapping must preserve positivity for quantities such as density and cloud water mixing ratio. 
  The cloud water mixing ratio is the amount of cloud water in air. 
  At each time step, the dynamics (advection) solutions, which are calculated on the advection mesh, are mapped to the reaction mesh to be used as input for the physics calculations. 
  The physics results are then mapped back to the dynamics to be used as input for the next time step. 
  Enforcing positivity alone may still lead to large oscillations and approximation errors.
  Using the DBI method will remove the large oscillations but will truncate any hidden extremum and may be too restrictive for high-order accuracy in some cases.
  For simulations in which different structured meshes are used and mapping is required, the errors from both the DBI and unconstrained PPI will propagate into other calculations and may even cause the simulation to fail.
  This section provides an analysis of the mapping error when interpolating from one mesh to another and back to the starting mesh.
  The mapping error is considered within time-dependent PDEs.
  For example, when interpolating the data values between the advection and reaction mesh in NEPTUNE, a mapping error is introduced in addition to the physics and time integration errors.
  The error in approximating a function $u(x)$ with the Newton polynomial $U_{n}(x)$ over the interval $I_{i}$ is 
  \begin{equation}\label{eq:interpError}
     E_{n}(x) = u(x)-U_{n}(x) = \frac{u^{(n+1)}(\xi)}{(n+1)!} \prod_{k=0}^{n}(x-x_{k}^{e}), \quad x \in I_{i}
  \end{equation}
  where $\xi \in [x_{n}^{l}, x_{n}^{r}]$.
  Given that $\xi$  and $u^{(n+1)}$ are not known, the local interpolation error in Equation (\ref{eq:interpError}) can approximated as follows:
  \begin{equation}\label{eq:error_approx}
    \tilde{E}_{n} = U[x_{n}^{l}\cdots x_{n}^{r}] \prod_{k=0}^{n}\Delta x_{k},
  \end{equation} 
   where
  \begin{equation*}\label{eq:Deltax}
    \Delta x_{k} = max \bigg( |x(i)-x_{k}^{e}|, |x(i+1)-x_{k}^{e}| \bigg).
  \end{equation*}
  The error approximation in Equation (\ref{eq:error_approx}) is based on the mean value theorem for divided differences, which states that there exist $\xi_{0} \in [x_{n}^{l}, x_{n}^{r}]$ such that 
  \begin{equation*}
    U[x_{n}^{l}\cdots x_{n}^{r}] =  \frac{u^{(n+1)}(\xi_{0})}{(n+1)!}.  
  \end{equation*}

  Equation (\ref{eq:error_approx}) approximates the local interpolation error for each interval when mapping from one set of points to another.
  To consider a mapping error for interpolating from one mesh to another and back to the starting mesh, let $\mathcal{M}_{A}$ and $\mathcal{M}_{R}$ be the advection and reaction mesh, respectively.
  In addition, let $I_{AR}$ and $I_{RA}$ be the interpolation operators that map a given set of data values from $\mathcal{M}_{A}$ to $\mathcal{M}_{R}$ and from $\mathcal{M}_{R}$  to $\mathcal{M}_{A}$, respectively. 
  We consider an advection-reaction problem where the advection part is calculated on $\mathcal{M}_{A}$ and the reaction on $\mathcal{M}_{R}$.
  A simple forward Euler time integration is used.  
  Let $\bar{u}_{\tau}$ and $\hat{u}_{\tau}$ be the approximate and the exact solution at time $\tau$.
  The dynamics/advection part is written as 
  \begin{equation}\label{eq:ti}
    \bar{u}^{1}_{\tau+\Delta \tau} = \bar{u}_{\tau} + \Delta \tau F\big( \bar{u}_{\tau} \big),
  \end{equation}
  and the physics/reaction $\bar{w}_{\tau+\Delta \tau}$ is expressed as
  \begin{equation}\label{eq:physics}
    \bar{w}_{\tau+\Delta \tau} = H\bar{u}^{1}_{\tau+\Delta \tau},
  \end{equation}
  where $H \bar{u}^{1}_{\tau+\Delta \tau}= I_{AR}G(I_{RA}\bar{u}^{1}_{\tau+\Delta \tau})$.
  Let $\bar{E}_{\tau+\Delta \tau}$ be the global space and time error accumulated up to $\tau+\Delta \tau$ after the advection in Equation (\ref{eq:ti}) and before mapping the solution values to $\mathcal{M}_{R}$.
  $\bar{E}_{\tau+\Delta \tau}$ does not include the mapping errors at $\tau+\Delta \tau$.
  The final solution after applying the operator $H$ in Equation (\ref{eq:physics}) is
  \begin{equation}\label{eq:u_bar}
    \bar{u}_{\tau+\Delta \tau} = \bar{u}^{1}_{\tau+\Delta \tau} + H\bar{u}^{1}_{\tau+\Delta \tau}.
  \end{equation}
  The true solution $\hat{u}_{\tau+ \Delta \tau}$ at the end of time step $\tau+\Delta \tau$ and after the mapping from $\mathcal{M}_{A}$ to $\mathcal{M}_{R}$ and back $\mathcal{M}_{A}$ to can be expressed as
  \begin{equation}\label{eq:u_hat}
    \hat{u}_{\tau+\Delta \tau} = \bar{u}^{1}_{\tau+ \Delta \tau} + \bar{E}_{\tau+\Delta \tau} +\hat{H} \big(\bar{u}^{1}_{\tau+\Delta \tau} + \bar{E}_{\tau+\Delta \tau} \big),
  \end{equation}
  where $\hat{H}$ is assumed to be the corresponding ``exact" operator for $H$.
  Subtracting Equation (\ref{eq:u_hat}) from (\ref{eq:u_bar}) gives an expression for the true error that can be written as
  \begin{equation*}
    E^{G}_{\tau+\Delta \tau} = \bar{E}_{\tau+\Delta \tau}+\hat{H} \big(\bar{u}^{1}_{\tau+\Delta \tau} + \bar{E}_{\tau+\Delta \tau} \big)- H\bar{u}^{1}_{\tau+\Delta \tau},
  \end{equation*} 
  where $E^{G}$ is the global space and time error including the mapping errors at $\tau+\Delta \tau$. 
  Adding and subtracting $H \big(\bar{u}^{1}_{\tau+\Delta \tau} + \bar{E}_{\tau+\Delta \tau} \big)$ yields
  \begin{equation*}
    E^{G}_{\tau+\Delta \tau} = \bar{E}_{\tau+ \Delta \tau}+\hat{H} \big(\bar{u}^{1}_{\tau+\Delta \tau} + \bar{E}_{\tau+\Delta \tau} \big)-H \big(\bar{u}^{1}_{\tau+\Delta \tau} + \bar{E}_{\tau+\Delta \tau} \big)
                     + H\big(\bar{u}^{1}_{\tau+\Delta \tau} + \bar{E}_{\tau+\Delta \tau} \big)-H\bar{u}^{1}_{\tau+\Delta \tau}.
  \end{equation*} 
  Using a Taylor expansion of $H\big(\bar{u}^{1}_{\tau+\Delta \tau} + \bar{E}_{\tau+\Delta \tau} \big)$ about $\bar{u}^{1}_{\tau+1}$ and dropping the high-order terms, we can approximate the total errors as 
  \begin{equation}\label{eq:total_error}
    E^{G}_{\tau+1} \approx \bar{E}_{\tau+\Delta \tau}+\hat{H} \big(\bar{u}^{1}_{\tau+\Delta \tau} + \bar{E}_{\tau+\Delta \tau} \big)-H \big(\bar{u}^{1}_{\tau+\Delta \tau} + \bar{E}_{\tau+\Delta \tau} \big)
                     + \frac{\partial H}{\partial u}(u^{1}_{\tau+\Delta \tau})E_{\tau+\Delta \tau}.
  \end{equation}
  The results in Equation (\ref{eq:total_error}) indicate that the total error is dependent on
  \begin{itemize}
    \item the existing global space and time error $\bar{E}_{\tau+\Delta \tau}$, which does not include the mapping error at $\tau+\Delta \tau$, 
    \item the mapping error $E^{M}_{\tau+\Delta \tau}$ at $\tau+\Delta \tau$,
          \begin{equation*}
            \begin{aligned}
            E^{M}_{\tau+\Delta \tau} &= \hat{H} \big(\bar{u}^{1}_{\tau+\Delta \tau} + \bar{E}_{\tau+\Delta \tau} \big)-H \big(\bar{u}^{1}_{\tau+\Delta \tau} + \bar{E}_{\tau+\Delta \tau} \big) \\
            &= \hat{H}\hat{u}^{1}_{\tau+\Delta \tau} - H \hat{u}^{1}_{\tau+\Delta \tau}, \textrm{ and }
            \end{aligned}
          \end{equation*}
    \item a multiplier of the existing global space and time error $\bar{E}_{\tau+\Delta \tau}$,
          \begin{equation*}
            E^{N}_{\tau+\Delta \tau} =\frac{\partial H}{\partial u}(u^{1}_{\tau+\Delta \tau})\bar{E}_{\tau+\Delta \tau}.
          \end{equation*}
  \end{itemize}
   Mapping data values from $\mathcal{M}_{A}$ to $\mathcal{M}_{R}$ and back to $\mathcal{M}_{A}$ introduces the interpolation errors that degrade the solution if $E^{M}_{\tau+\Delta \tau}$ is greater than the existing global space and time error $\bar{E}_{\tau+ \Delta \tau}$. 
   This problem is resolved when the mapping error is kept smaller than the existing global space and time error.
   Similar ideas in the context of time dependent differential equations are explored in \cite{doi:10.1137/050622092, doi:10.1137/0912031, BERZINS1998117}. 
   The studies in \cite{doi:10.1137/050622092} and \cite{doi:10.1137/0912031} develop strategies for balancing the space and time error for better error control and improved performance while \cite{BERZINS1998117} show that in mesh adaptivity the spatial interpolation error must be controlled and kept smaller than the temporal error.

    \subsection{1D Modified Runge Function}
    \label{subsec:runge}
      This example is based on the modified version of the Runge function defined in Equation (\ref{eq:runge-1d}) and two meshes that emulate the dynamics/advection ($\mathcal{M}_{A}$) and physics/reaction ($\mathcal{M}_{R}$) meshes used in NEPTUNE.
      No actual PDE is solved in this example.
      Here, we consider the advection time step of the trivial case where the identity operator is used to represent the advection and reaction.
      The function $f_{1}(x)$ is evaluated on advection mesh $\mathcal{M}_{A}$ to create the initial data values.
      Given that the identity operator is used for both the advection and reaction, these initial data values are mapped to the reaction mesh $\mathcal{M}_{R}$ and back to the starting mesh $M_{A}$.
      Using the identity operator allows for a study of the mapping error without the influence of the advection and reaction. 

      Table \ref{tab:f4-map} shows $L^{2}$-norms of the mapping errors over the grid points for $f_{1}(x)$ when using the PCHIP, MQSI, DBI, and PPI methods to map the data values from the advection mesh to the reaction mesh and back to the advection mesh.
      For $N=64$ points, increasing the interpolant degree does not significantly improve the approximation.
      The global error is dominated by the local error in the regions with steep gradients that are to the left and right of the peak at $x=0$.
      The mapping errors can be improved by increasing the resolution and adding more points in the regions with steep gradients.
      The resolution is increased by adding one or three uniformly spaced points in each interval from the initial profile with 64 points.
      Increasing the resolution leads to better approximations when mapping data values between both meshes, and the error decreases as we increase the polynomial degree from $3$ to $7$.
      This example demonstrates that in cases with steep gradients, using the PPI method with high-order interpolants may not significantly improve the approximation unless there is sufficient resolution.
      In order to benefit from the positivity and the high-order interpolants, it is important to be in the regime where the problem has sufficiently many points to observe convergence as the polynomial degree increases. 
      Overall, the PPI method leads to smaller errors compared to the other methods as the resolution and polynomial degree increase.
%      %
    \begin{table}[H]
        \centering
        \begin{tabular}{ c c c c c c c c c c c c c}
          \hline
          \hline
          $N$  && PCHIP   &&  MQSI  && \multicolumn{3}{c}{DBI}            &&   \multicolumn{3}{c}{PPI}         \\
               && $\mathcal{P}_{3}$ && $\mathcal{P}_{5}$ && $\mathcal{P}_{3}$ & $\mathcal{P}_{5}$ & $\mathcal{P}_{7}$                                        
                                    && $\mathcal{P}_{3}$ & $\mathcal{P}_{5}$ & $\mathcal{P}_{7}$ \\
          \hline
          64       && 2.92E-03  &&  1.93E-03  &&  4.93E-03  &  4.78E-04  &  7.12E-05  &&  3.99E-03  &  3.13E-04  &  2.85E-05   \\
          127      && 3.81E-04  &&  5.57E-04  &&  3.58E-03  &  3.81E-04  &  6.50E-05  &&  2.85E-03  &  1.02E-04  &  3.65E-06   \\
          253      && 6.71E-05  &&  1.48E-04  &&  3.41E-03  &  3.69E-04  &  6.49E-05  &&  2.46E-03  &  3.50E-05  &  8.62E-07   \\
          \hline
          \hline
        \end{tabular}                                                                                  
        \caption{$L^{2}$-norm of mapping errors for the modified Runge function $f_{4}(x)$ when using the PCHIP, MQSI, DBI, and PPI methods to map the data values from the advection mesh to the reaction mesh and back to the advection mesh.
                  The target polynomials are set to $d=3$, $d=5$, and $d=7$.
                 $N$ represents the number of input points used for both meshes.
                 The parameter $st$ is set to 3.}
        \label{tab:f4-map}
      \end{table}

    \subsection{TWP-ICE Example}
    \label{subsec:twp-ice}
      This study uses the tropical warm pool international cloud experiment (TWP-ICE) test case from the common community physics package (CCPP) \cite{ccpp-v4}.
      The input mesh for the simulation is configured to emulate a vertical column in NEPTUNE.
      The simulation result at time $t=1440$ sec is extracted and scaled, and used to evaluate different interpolation approaches when mapping solution values between advection and reaction meshes. 
      The domain and range are scaled to $[-1,1]$ and $[0,1]$, respectively. 
      This study considers the cloud water mixing ratio profile, which represents the amount of cloud water in air. 
      The extracted profile is then fitted using a radial basis function interpolation to construct an analytical function that can be used as the starting point of the mapping evaluation.
      The radial basis function is based on multiquadrics:
      \begin{equation*}
          b_{i} = \sqrt{1 + (\epsilon|x-x_{i}|)^{2}}.
      \end{equation*}
      The parameter $\epsilon$ is approximated using cross validation ~\cite{fasshauer2007choosing}.
      The initial values are obtained by evaluating the analytical function on the advection mesh.
      These values are then mapped to the reaction mesh and back to the advection mesh.   

      Table \ref{tab:twp-ice-map} shows $L^{2}-$norms of the mapping errors for the extracted profile when using the PCHIP, DBI, and PPI methods to map the data values from the advection to physics mesh and back to advection mesh.
      For $N=64$, the global error is dominated by the local error at a few points located in the regions with steep gradients.
      Increasing the polynomial degree does not significantly improve the approximation compared to using PCHIP for $N=64$.
      More points are required to better approximate the underlying profile in the regions with steep gradients.     
      The resolutions are increased by adding one and three uniformly spaced points in each interval from the initial $N=64$ mesh points.
      Table \ref{tab:twp-ice-map} shows that with the increased resolution, the approximation improves as the polynomial degree increases.
      The number of points used in each region with steep gradients increased as more points were added.  
      This example provides an application example using simulation data from TWP-ICE.
      In cases of coarse resolution (64 points), the PCHIP and MQSI results are comparable and going to higher degree interpolants doesn't significantly improve the approximation for DBI and PPI.
      The approximation improves with higher degree interpolants when the resolution is increased, as shown in Table \ref{tab:twp-ice-map}.
      The results from this experiment suggest that increasing the resolution is needed for the mapping between meshes to benefit from the high-order interpolants from the PPI methods.
      \begin{table}[H]
        \centering
        \begin{tabular}{ c c c c c c c c c c c c c}
          \hline
          \hline
          $N$  && PCHIP   &&  MQSI  && \multicolumn{3}{c}{DBI}            &&   \multicolumn{3}{c}{PPI}         \\
               && $\mathcal{P}_{3}$ && $\mathcal{P}_{5}$ && $\mathcal{P}_{3}$ & $\mathcal{P}_{5}$ & $\mathcal{P}_{7}$                         
                                    && $\mathcal{P}_{3}$ & $\mathcal{P}_{5}$ & $\mathcal{P}_{7}$        \\
         \hline
         64       && 4.66E-03  &&  3.25E-03  &&  1.17E-02  &  3.15E-03  &  6.77E-04  &&  5.11E-03  &  9.86E-04  &  1.12E-04   \\
         127      && 1.56E-03  &&  2.05E-03  &&  1.12E-02  &  3.10E-03  &  6.49E-04  &&  2.30E-03  &  3.10E-04  &  1.83E-05   \\
         253      && 4.89E-04  &&  6.92E-04  &&  1.11E-02  &  3.11E-03  &  6.46E-04  &&  1.41E-03  &  1.45E-04  &  3.70E-06   \\
          \hline
          \hline
        \end{tabular}                                                                                  
        \caption{$L^{2}-$norm of mapping errors for the TWP-ICE profile when using the PCHIP, MQSI, DBI, and PPI methods to map the data values from the advection mesh to the reaction mesh and back to the advection mesh.
                  The target polynomials are set to $d=3$, $d=5$, and $d=7$.
                 $N$ represents the number of input points used for both meshes.
                 The parameter $st$ is set to 3.}
        \label{tab:twp-ice-map}
      \end{table}

 \subsection{BOMEX Example}
      Here, a maritime shallow convection example based on the 1D Barbados Oceanographic and Meteorological Experiment (BOMEX) \cite{friedman1970essa} is used to study the effects of the different interpolation methods in an application.
      BOMEX is a single-column shallow convection test case used to measure and study changes in temperature,  wind, water vapor mixing ratio, rain water mixing ratio, and cloud water mixing ratio.
      The mixing ratios represent the amount of water vapor, rain water, and cloud water in air. %, the properties of heat, moisture, and momentum.
      In this simulation, the dynamics/advection is modeled by 
      \begin{equation*}\label{eq:bomex}
          \frac{\partial X}{\partial t} = \mathcal{L}X,
      \end{equation*}
      where $X$ is the state variable that contains the wind, temperature, water vapor mixing ratio, cloud water mixing ratio, and rain water mixing ratio.
      The dynamics are approximated using fifth-order weighted essentially nonoscillatory (WENO) and third-order Runge-Kutta methods ~\cite{doi:10.1080/1061856031000104851}.
      The physics component of the simulation uses the hybrid eddy-diffusivity mass-flux and free atmospheric turbulence (hybrid EDMF) and Kessler microphysics schemes from \cite{ccpp-v4} to alter the results of the dynamics and incorporate the physics properties.
      The dynamics and physics results are calculated on different meshes.
      At each time step, the dynamics are calculated on the advection mesh $\mathcal{M}_{A}$, and the results are interpolated to the reaction mesh $\mathcal{M}_{R}$ for the use of the physics routines.
      The physics terms are calculated using the reaction mesh $\mathcal{M}_{R}$, and the results are interpolated back to the advection mesh $\mathcal{M}_{A}$.

\begin{figure}[H]
         \begin{subfigure}{0.4\textwidth}
         \centering
         \includegraphics[scale = 0.10]{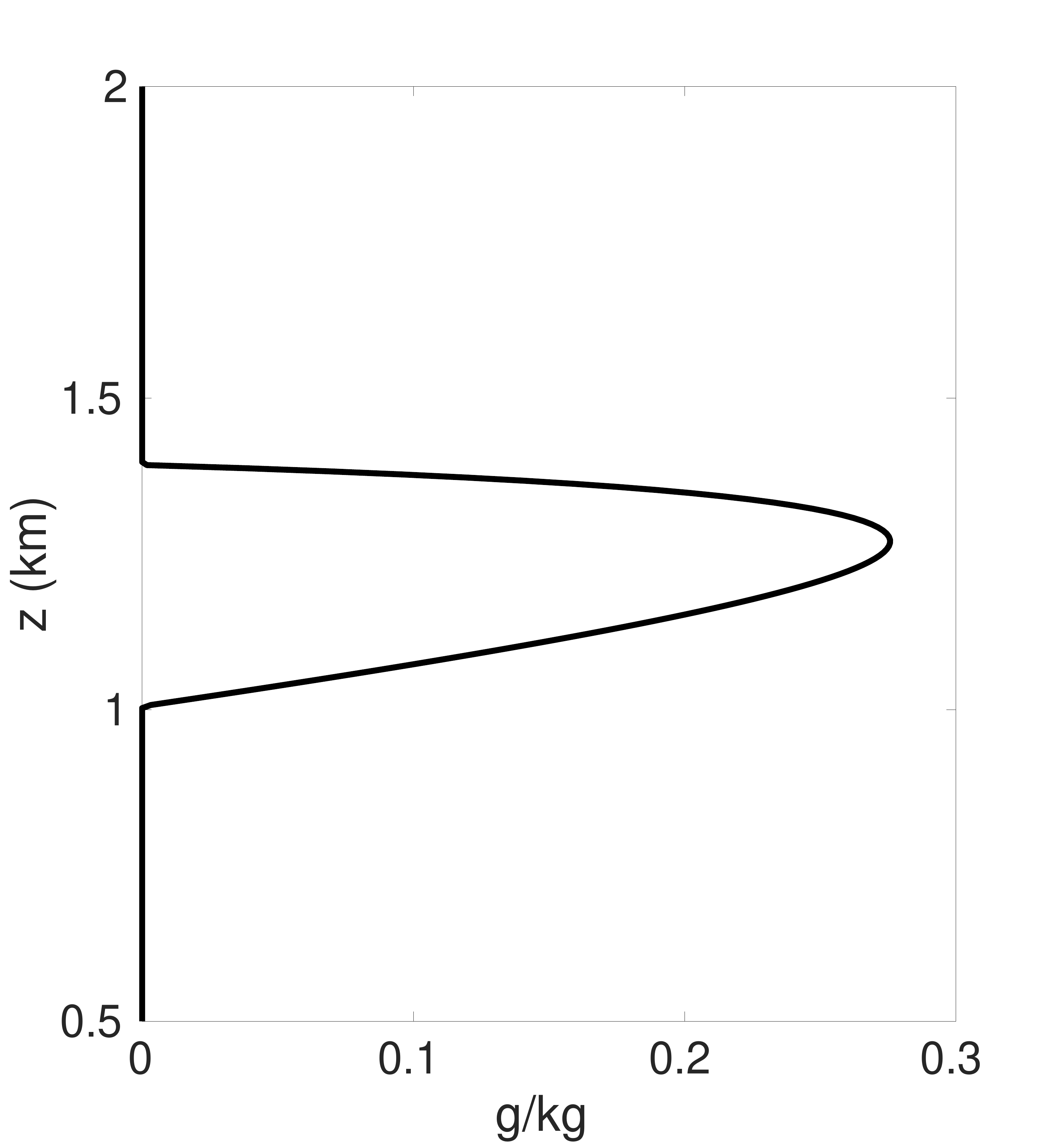}
         \subcaption{Target (no mapping required).}
         \label{subfig:qcweno_physics}
         \end{subfigure}
         \begin{subfigure}{0.4\textwidth}
         \centering
         \includegraphics[scale = 0.10]{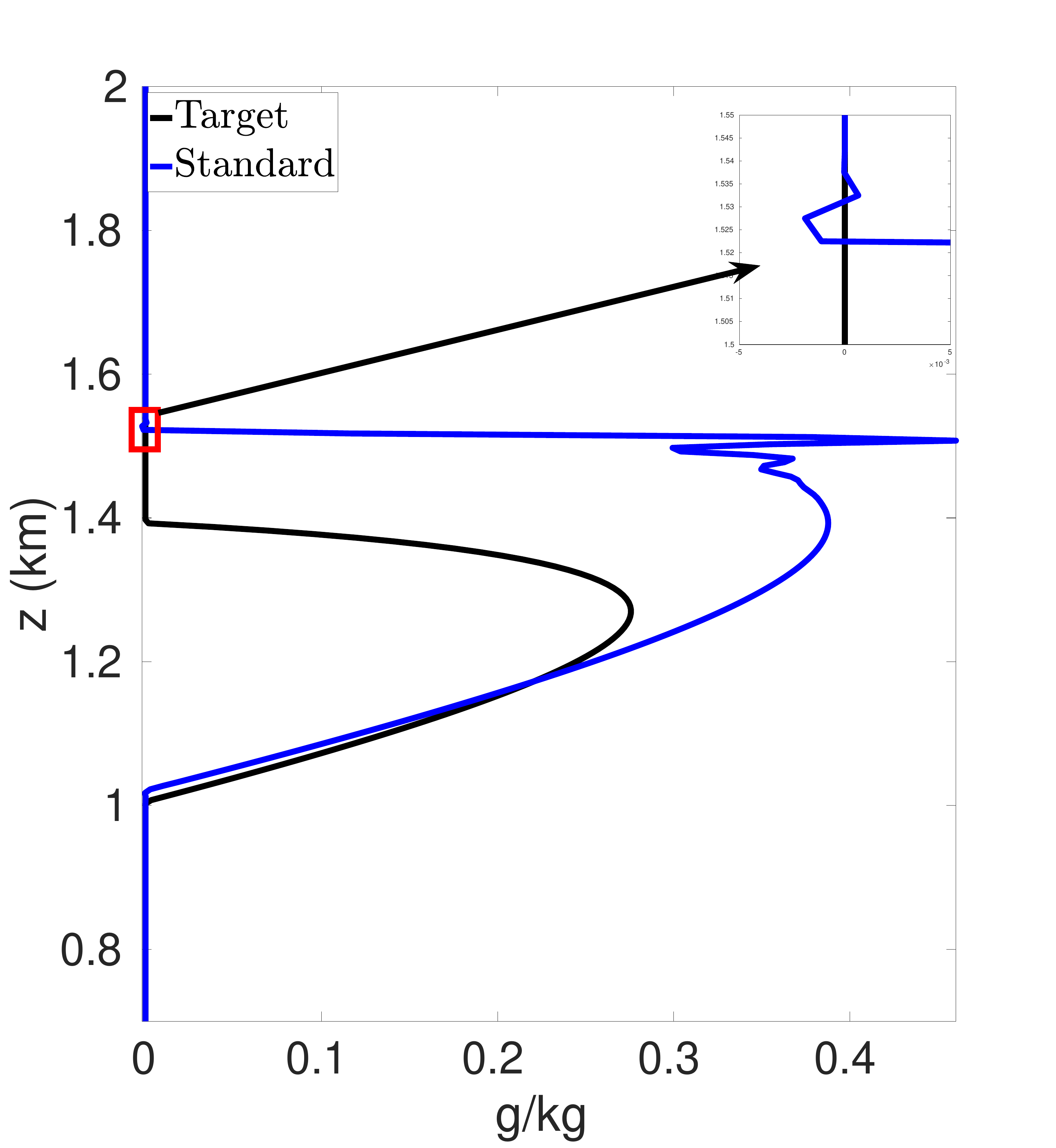}
         \subcaption{Standard interpolation.}
         \label{subfig:qcweno_std}
         \end{subfigure}
         \begin{subfigure}{0.4\textwidth}
         \centering
         \includegraphics[scale = 0.10]{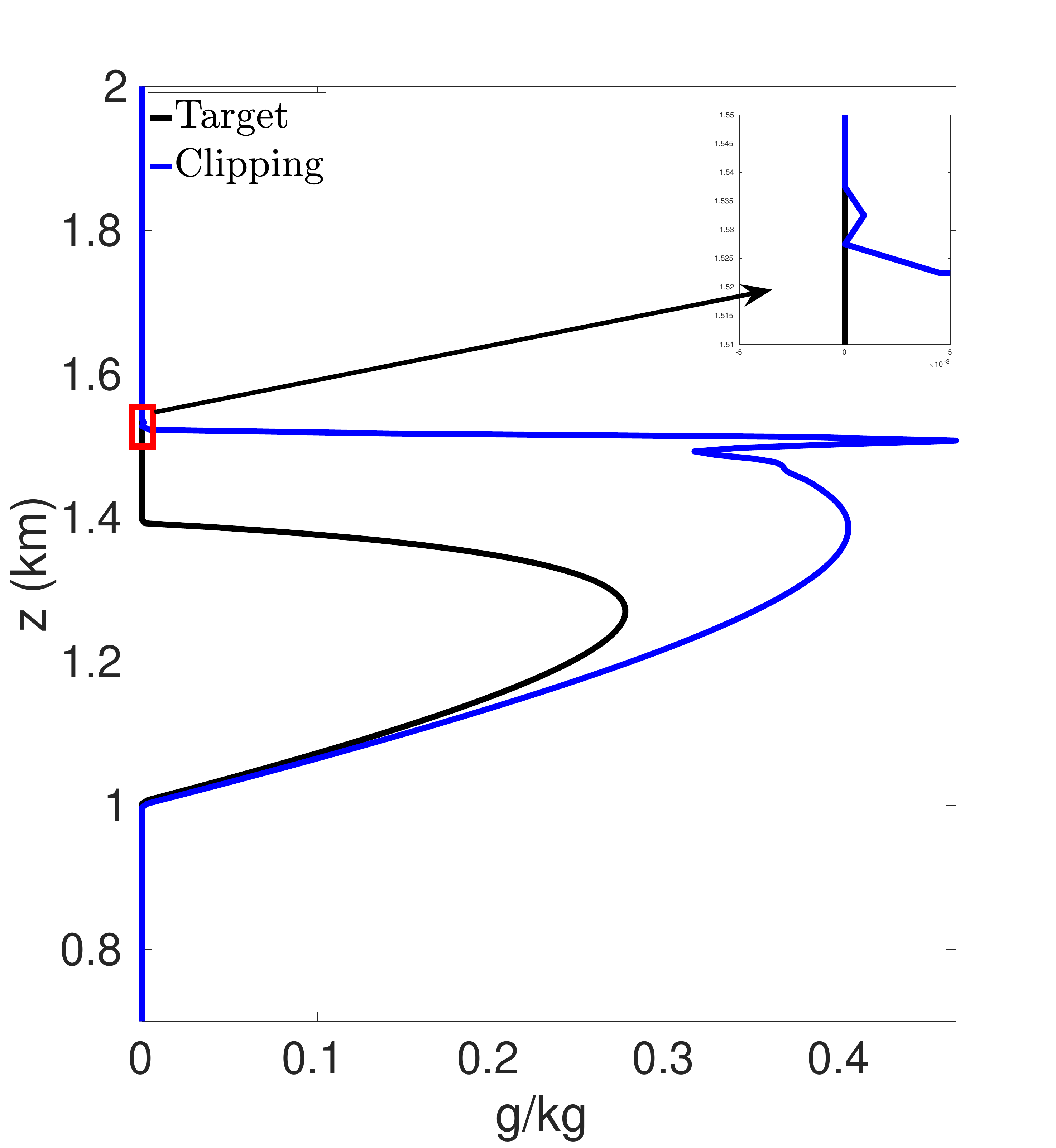}
         \subcaption{Standard interpolation with ``clipping".}
         \label{subfig:qcweno_clip}
         \end{subfigure}
         \begin{subfigure}{0.4\textwidth}
         \centering
         \includegraphics[scale = 0.10]{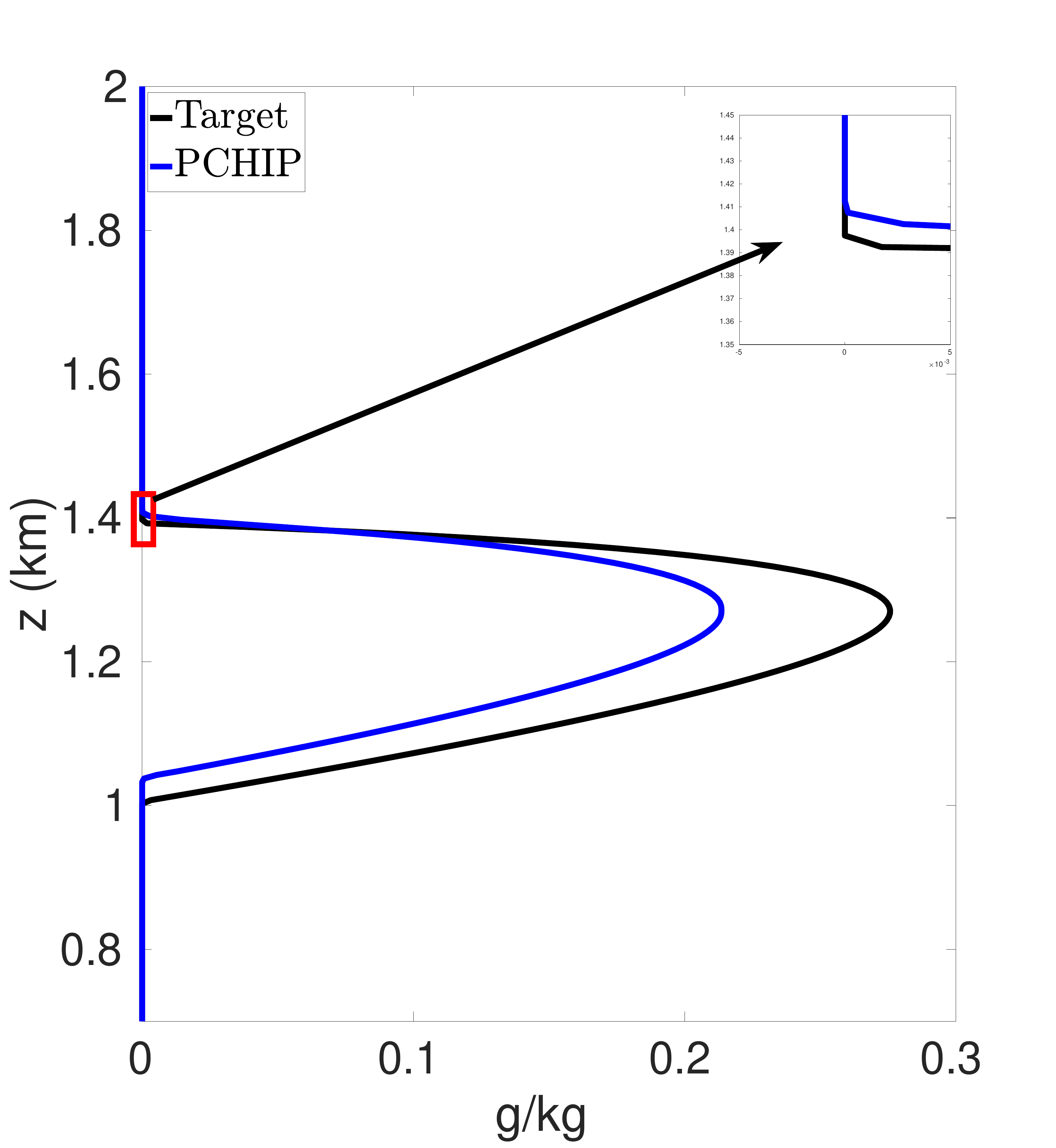}
         \subcaption{PCHIP.}
         \label{subfig:qcweno_pchip}
         \end{subfigure}
         \caption{Cloud mixing ratio $q_{c}$ profile from the BOMEX test case at $t= 5$ hours with $nz= 600$ points.
                 A fifth-order WENO  and third-order Runge-Kutta schemes with $CFL= 0.1$ are used for the dynamics (advection).
                 \ref{subfig:qcweno_physics}.
                 The black plot in \ref{subfig:qcweno_std}, \ref{subfig:qcweno_clip}, and \ref{subfig:qcweno_pchip} represents the target profile where the same mesh is used for the dynamics and physics calculations. 
                 In \ref{subfig:qcweno_std}, \ref{subfig:qcweno_clip}, and  \ref{subfig:qcweno_pchip}, the profiles in blue use different meshes for the dynamics and physics calculations which require mapping the solution values between both meshes. 
                 A standard polynomial interpolation, a standard polynomial interpolation with ``clipping", and PCHIP methods are used for the mapping in \ref{subfig:qcweno_std}, \ref{subfig:qcweno_clip}, and \ref{subfig:qcweno_pchip}, respectively.} 
         \label{fig:qcweno_std_pchip}
      \end{figure}

\begin{figure}[H]
        \vspace{-3mm}
         \begin{subfigure}{0.4\textwidth}
         \centering
         \includegraphics[scale = 0.10]{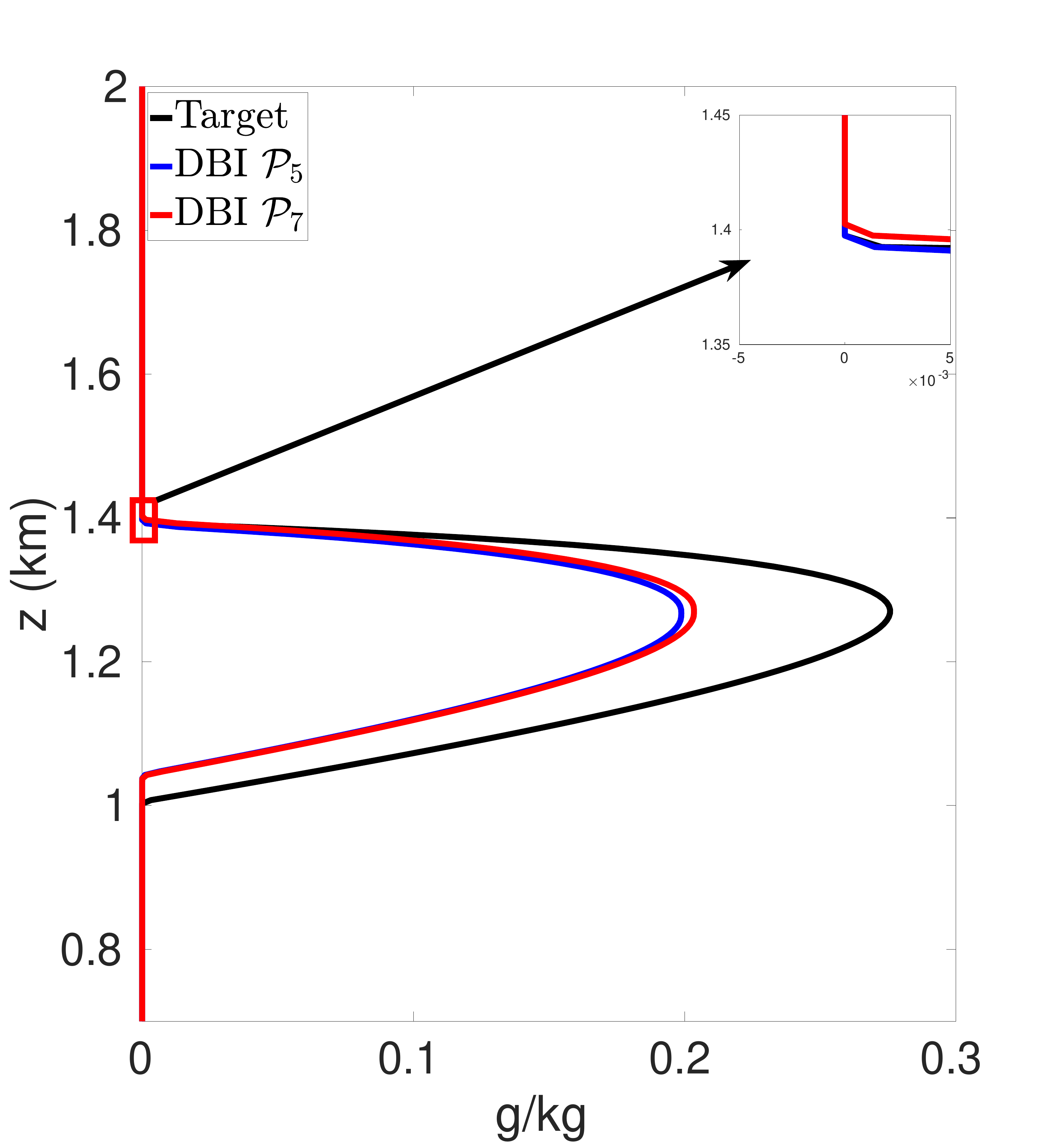}
         \subcaption{}
         \label{subfig:qcweno_dbi_s1}
         \end{subfigure}
         \begin{subfigure}{0.4\textwidth}
         \centering
         \includegraphics[scale = 0.10]{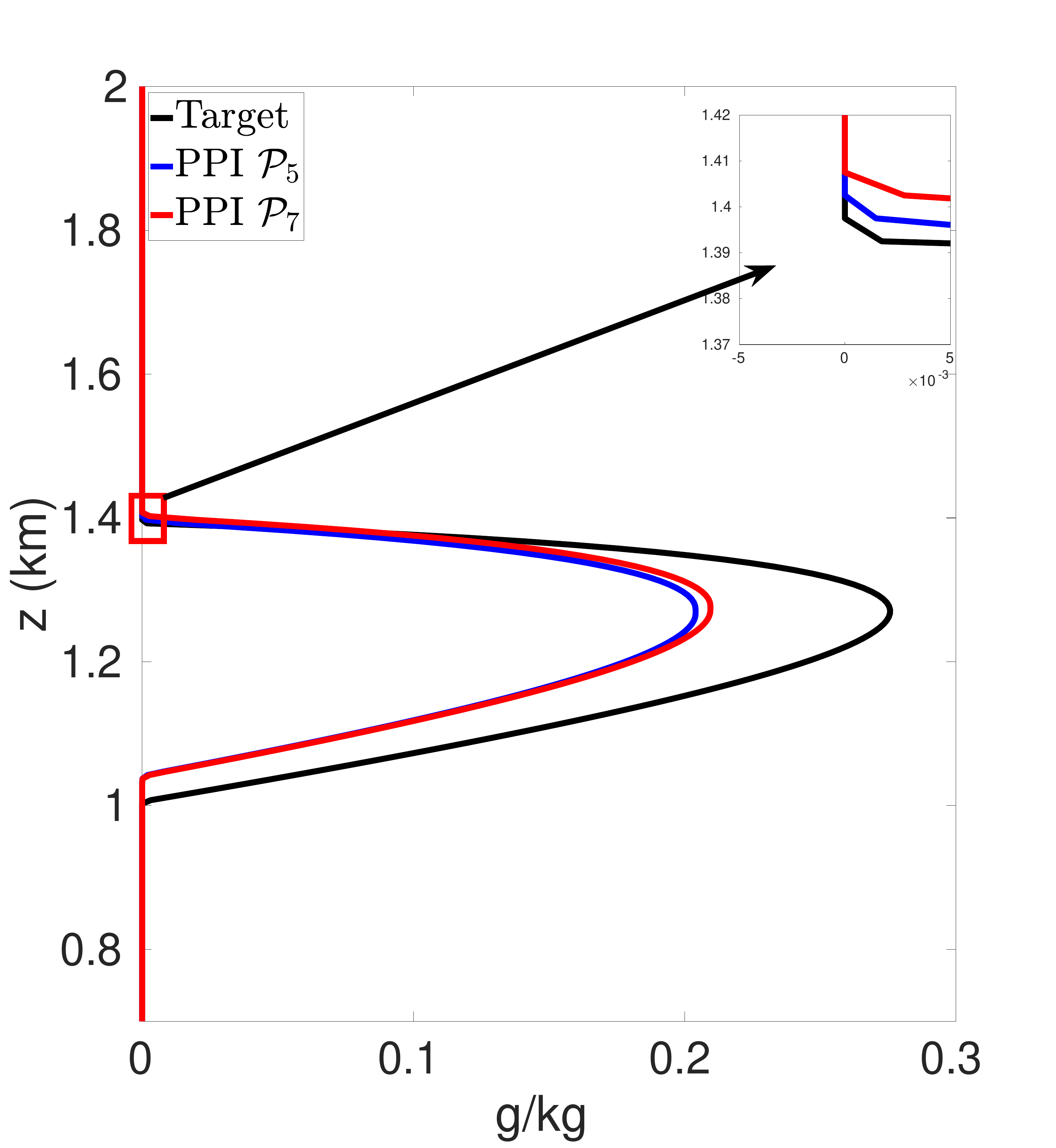}
         \subcaption{}
         \label{subfig:qcweno_ppi_s1}
         \end{subfigure}
         \vspace{-0.5mm}
         \begin{subfigure}{0.4\textwidth}
         \centering
         \includegraphics[scale = 0.10]{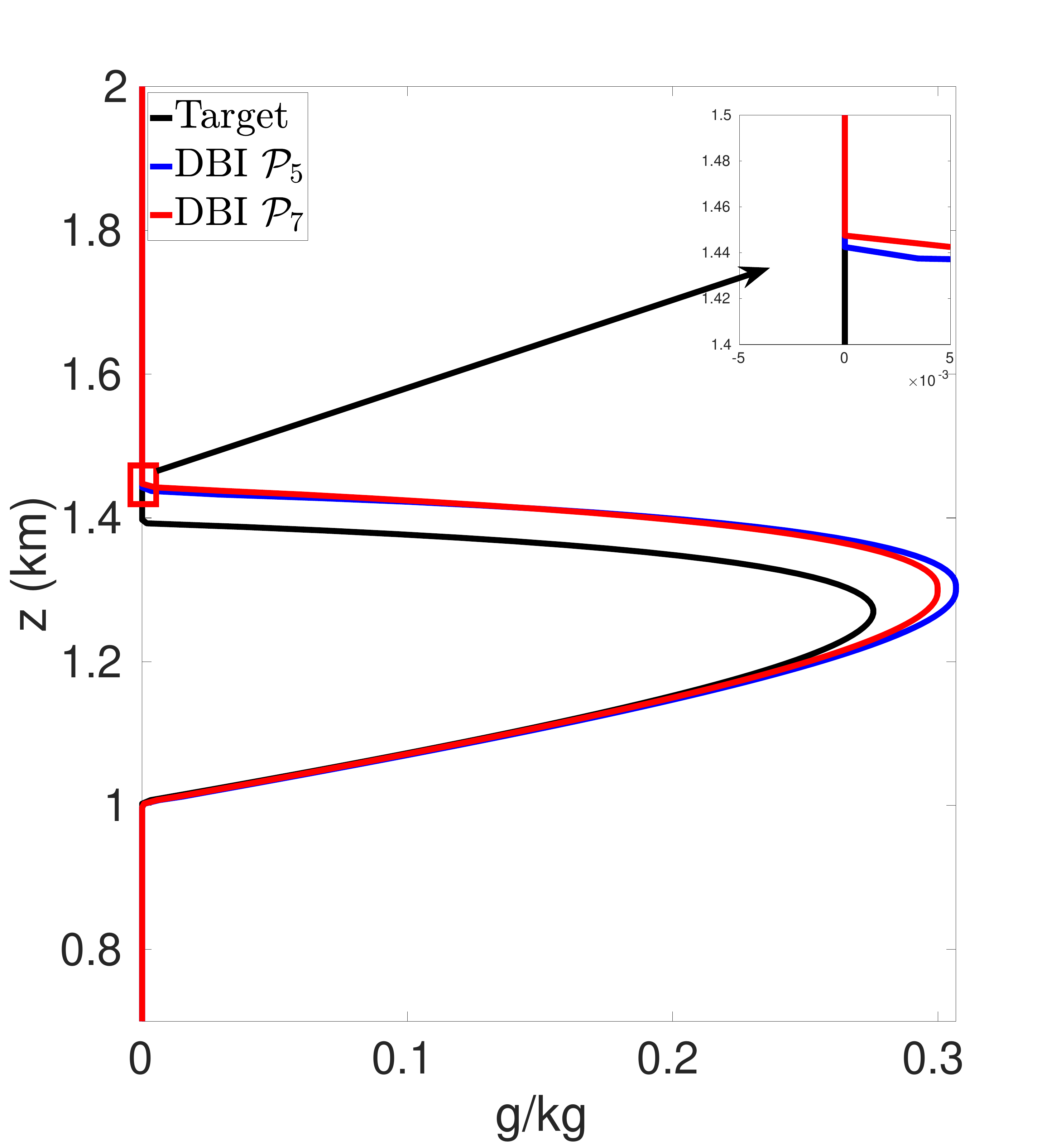}
         \subcaption{}
         \label{subfig:qcweno_dbi_s2}
         \end{subfigure}
         \begin{subfigure}{0.4\textwidth}
         \centering
         \includegraphics[scale = 0.10]{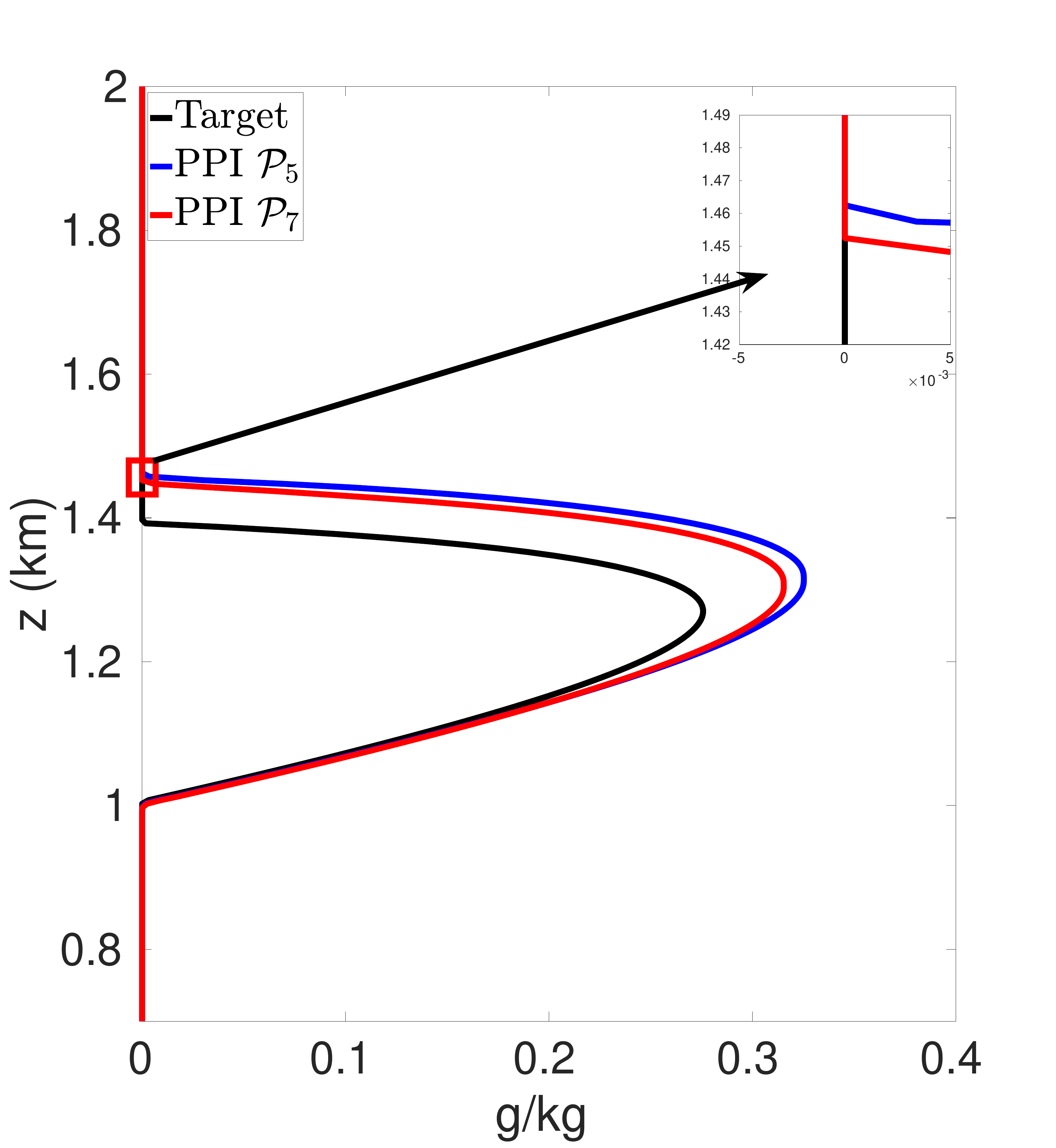}
         \subcaption{}
         \label{subfig:qcweno_ppi_s2}
         \end{subfigure}
         \vspace{-0.5mm}
         \begin{subfigure}{0.4\textwidth}
         \centering
         \includegraphics[scale = 0.10]{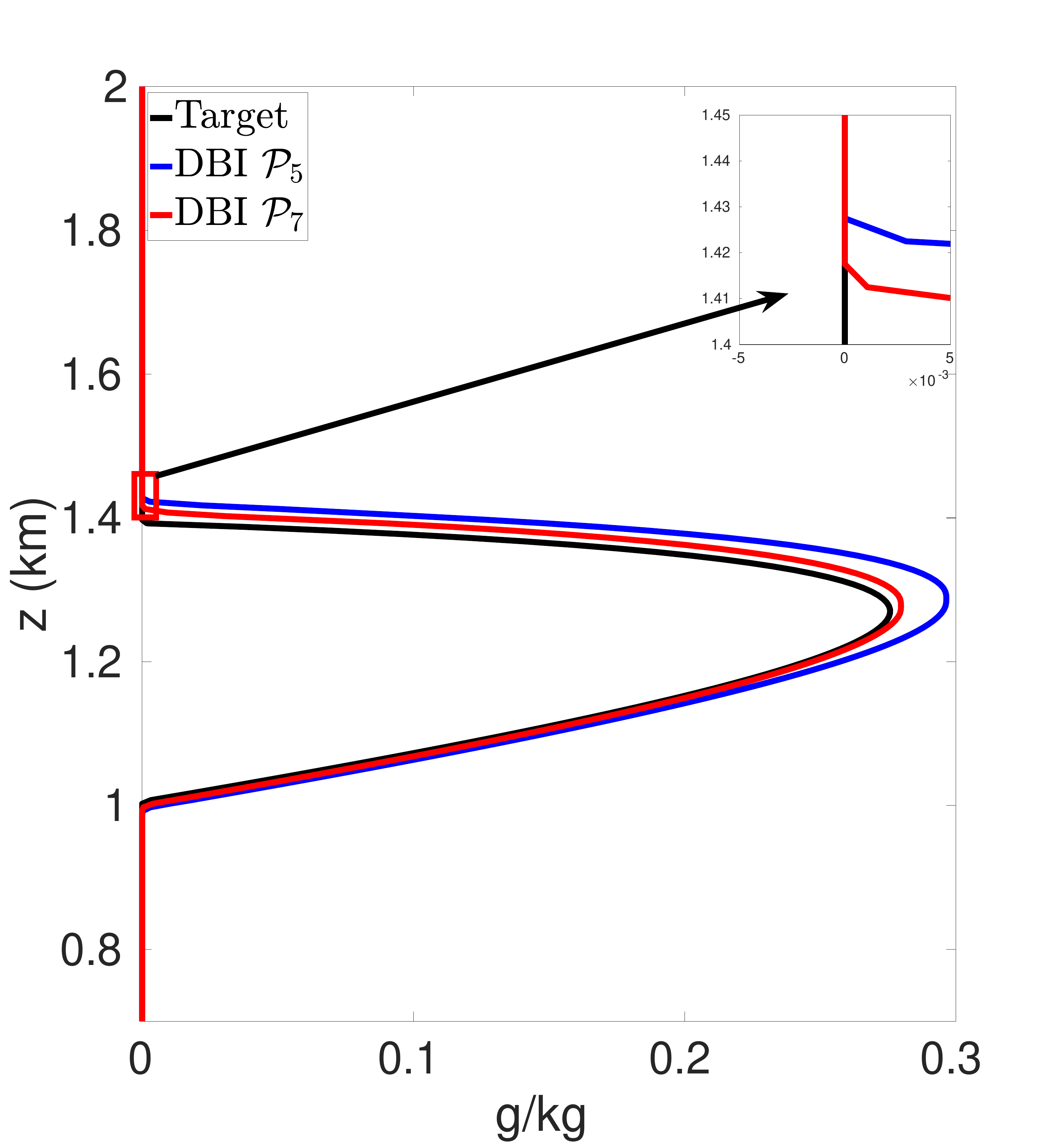}
         \subcaption{}
         \label{subfig:qcweno_dbi_s3}
         \end{subfigure}
         \begin{subfigure}{0.4\textwidth}
         \centering
         \includegraphics[scale = 0.10]{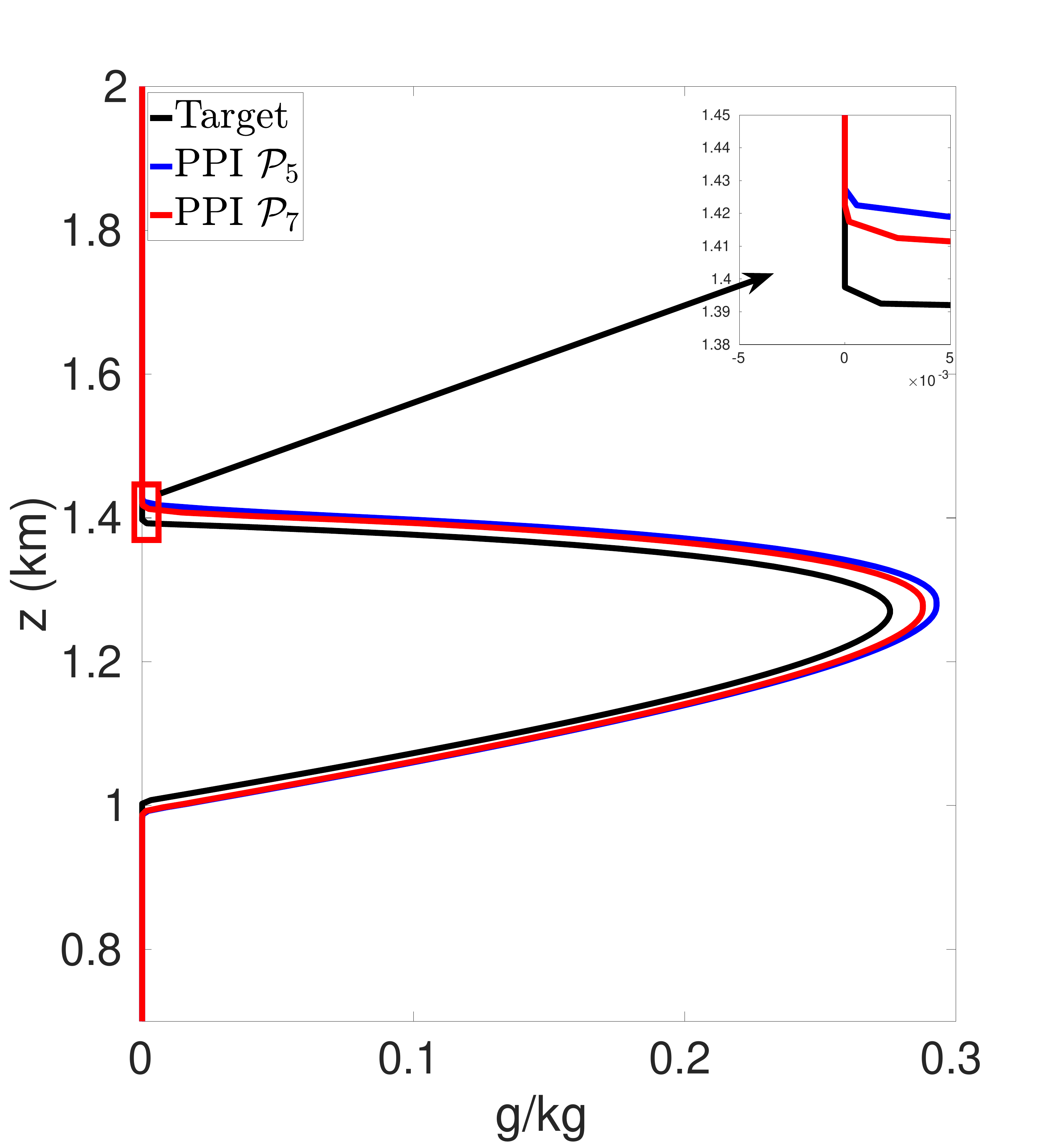}
         \subcaption{}
         \label{subfig:qcweno_ppi_s3}
         \end{subfigure}
         \vspace{-1mm}
         \caption{Cloud mixing ratio $q_{c}$ profile from the BOMEX test case at $t= 5$ hours with $nz= 600$ points with $\epsilon_{0}=\epsilon_{1}=10^{-5}$. 
                  The profile in black is the target solution. 
                  The profiles on the left (\ref{subfig:qcweno_dbi_s1}, \ref{subfig:qcweno_dbi_s2}, \ref{subfig:qcweno_dbi_s3}) and right (\ref{subfig:qcweno_ppi_s1}, \ref{subfig:qcweno_ppi_s2}, \ref{subfig:qcweno_ppi_s3}) are obtained using the DBI and PPI methods, respectively, to map solution values between meshes.
                  The maximum polynomial degrees are set to $d=5$ and $d=7$ for the blue and red plots, respectively.
                  The input parameter $st$ is set to $1$, $2$, and $3$ for the advection (\ref{subfig:qcweno_dbi_s1}, \ref{subfig:qcweno_ppi_s1}), reaction (\ref{subfig:qcweno_dbi_s2}, \ref{subfig:qcweno_ppi_s2}), and third (\ref{subfig:qcweno_dbi_s3}, \ref{subfig:qcweno_dbi_s3}, \ref{subfig:qcweno_ppi_s3}) row, respectively. 
                  A fifth-order WENO  and third-order Runge-Kutta schemes with $CFL= 0.1$ are used for the dynamics (advection).}
         \label{fig:qcweno_dbi_ppi}
      \end{figure}

      As in \cite{Rotstayn2000}, let $q_{c}$ be the cloud water mixing ratio profile in the different experiments.
      Figures \ref{subfig:qcweno_physics} - \ref{subfig:qcweno_ppi_s3} show the cloud mixing ratio profile  $q_{c}$ at $t=5$ hours that is used as input for the physics routines.
      The physics calculations require positive input values for $q_{c}$.
      Figure \ref{subfig:qcweno_physics} shows the target profile for $q_{c}$.
      This target profile is obtained by using the same mesh for both dynamics and physics calculations where mapping is not required and $q_{c}$ remains positive during the simulation. 
      In addition, as the temporal and spatial resolution increases, $q_{c}$ converges to the profile shown in Figure \ref{subfig:qcweno_physics}.
      Figures \ref{subfig:qcweno_std} - \ref{subfig:qcweno_ppi_s3} are used to investigate different interpolation methods for mapping the solution values between meshes in the case where the dynamics and physics are calculated using different meshes. 

      Figure \ref{subfig:qcweno_std} shows the cloud mixing ratio profiles $q_{c}$ for the target and approximated solution at $t= 5$ hours.
      In the case of the approximated solution, a fifth-order standard polynomial interpolation is used when mapping between the advection and reaction meshes.
      For a given interval $I_{i}$, the polynomial interpolant is constructed using the stencil $\mathcal{V}_{4}=\{x_{i-2}, x_{i-1}, x_{i}, x_{i+1}, x_{i+2}, x_{i+3}\}$.
      At the boundary and nearby boundary intervals, the stencil $\mathcal{V}_{4}$ is biased toward the interior of the domain. 
      The results in Figure \ref{subfig:qcweno_std} demonstrate that using the standard polynomial interpolation leads to oscillations, negative values, and an overestimation of the peak and total cloud mixing ratio of the profile $q_{c}$.
      Using standard polynomial interpolation leads to an overproduction of the total cloud mixing ratio by $93.45\%$.
      The peak is $max(q_{c}) = 0.46 g/kg$, which is larger than the target peak $max(q_{c})= 0.28 g/kg$. 

      The negative values in Figure \ref{subfig:qcweno_std} can be removed via ``clipping", which is a procedure that consists of removing the negative values by setting them to zero. % \cite{skamrock}.
      Figure \ref{subfig:qcweno_clip} shows the cloud mixing ratio profiles for the target solution and an approximated solution that uses ``clipping" to remove the negative values at each time step.
      The approximated solution uses a standard interpolation to map the data values from one mesh to another. 
      The interpolant for each interval is constructed using the stencil $\mathcal{V}_{4}=\{x_{i-2}, x_{i-1}, x_{i}, x_{i+1}, x_{i+2}, x_{i+3}\}$ with a fifth-order polynomial.
      Once the interpolation is completed, ``clipping" is used to remove the negative values.
      Figure \ref{subfig:qcweno_clip} shows that using ``clipping"  still allows for oscillations and a positive bias in the prediction of the cloud mixing ratio $q_{c}$. 
      The total cloud mixing ratio is $2.09$ times greater than the target solution, and the peak $max(q_{c})=0.46 g/kg$ is larger than the target peak $max(q_{c})=0.28 g/kg$.

      Using PCHIP to map between the advection and reaction meshes eliminates the negative values, removes oscillations, and reduces the positive bias in the cloud mixing ratio prediction compared to the standard interpolation with and without ``clipping".
      Figure \ref{subfig:qcweno_pchip} shows the target profile $q_{c}$ and an approximated profile that uses PCHIP for mapping solution values between advection and reaction meshes.
      %
      %Using PCHIP reduces the positive bias prediction but does not remove it.
      %
      The total cloud mixing ratio is now $27.21\%$ less than the target with a peak $max(q_{c})= 0.21 g/kg$.
      In the BOMEX test case, NEPTUNE, and similar codes, using PCHIP for mapping data values from one mesh to another can degrade the high-order accuracy obtained from the high-order methods used for the dynamics calculations.
      PCHIP is only third-order whereas the dynamics calculations use a fifth-order method.
      This limitation can be addressed via high-order DBI and PPI. 
      Here, the MQSI method is not used because it oversmoothes the state variable $X$ at each time step and leads to no production of cloud mixing ratio.
  
      Figures \ref{subfig:qcweno_dbi_s1}-\ref{subfig:qcweno_ppi_s3} show cloud mixing ratio profiles for the target and approximated solutions that use the DBI and PPI methods to map the solution values between meshes.
      The maximum polynomial degree for the DBI and PPI methods is set to $5$ and $7$, and the parameters $\epsilon_{0}$ and $\epsilon_{1}$ are both set a value of $10^{-5}$.
      For larger values of $\epsilon_{0}$ and $\epsilon_{1}$, the PPI approach introduces oscillations that lead to positive bias prediction of the cloud mixing ratio.
      These oscillations are caused by the relaxed nature of the PPI approach, which still allows the interpolants to oscillate while remaining positive.
      The positive bias and oscillations can be removed using the DBI or PPI method with small values for $\epsilon_{0}$ and $\epsilon_{1}$.
      When using the PPI method for mapping, the total amount of the cloud mixing ratio is less than the target for $st=1$ and more than the target for $st=2$ and $st=3$.
      Figures \ref{subfig:qcweno_dbi_s1}-\ref{subfig:qcweno_ppi_s3} show that using the  DBI and PPI methods with $\epsilon_{0}=\epsilon_{1}=10^{-5}$ to map data values between the advection and reaction meshes eliminates the negative values, removes the oscillations, and significantly reduces the positive bias in the cloud mixing ratio prediction.
      Using the DBI and PPI methods leads to a better approximation of the peak value of the total cloud mixing ratio compared to using the standard interpolation and  PCHIP approaches. 
      The best approximation of the total amount of the cloud mixing ratio is with the DBI method, which is $7.57\%$ more than the target with a peak of $max(q_{c}) = 0.28  g/kg$.

      In summary, using DBI and PPI methods to map data values between both the advection and reaction meshes produces better approximation results compared to the standard interpolation and  PCHIP methods.
      Tables \ref{tab:qc0} and \ref{tab:qc1} provide a summary of the maximum values and the total amount of cloud mixing ratios for each case.
      The DBI and PPI methods with a target polynomial set to $d=7$ lead to a better approximation of the peak and the total cloud mixing ratios compared to the standard interpolation and PCHIP approaches.
      The results from Tables \ref{tab:qc0} and \ref{tab:qc1} indicate that the DBI method is the most suitable approach to map data values between meshes for the BOMEX test case.
      This study provided an example demonstrating how to use the DBI and PPI methods for mapping data values between meshes in the context of NWP.
      The BOMEX example also demonstrated that positivity alone may not be sufficient to remove oscillations in the solution, and the interpolants may need to be constrained to be between the data values for a better approximation.
 
    \begin{table}[H]
      \centering
      \begin{tabular}{ c c c c c}
        \hline
        \hline
                        & Target    & STD     & Clipping & PCHIP   \\
        \hline
        maximum $q_{c}$ & 0.28      & 0.46    & 0.46     & 0.21    \\
        total $q_{c}$   & 69.82   & 135.07   & 145.89   & 50.82    \\
        \hline
        \hline
      \end{tabular}                                                                                  
      \caption{Maximum values of $q_{c}$ and the total amount of the cloud mixing ratio at $t = 5$ hours with $nz=600$ points.
               The total amount of the cloud mixing ratio is calculated by estimating the integral $q_{c}$.
               The units of $q_{c}$ are $g/kg$.}
      \label{tab:qc0}
    \end{table}
    \begin{table}[H]
      \centering
      \begin{tabular}{ c c c c c c c c c c}
        \hline
        \hline
                          & $st=1$  &  $st=2$ & $st=3$  && $st=1$ &  $st=2$ & $st=3$  &&\\
                          & \multicolumn{3}{c}{$\mathcal{P}_{5}$} && \multicolumn{3}{c}{$\mathcal{P}_{7}$} && Target\\
        \hline
                          & \multicolumn{7}{c}{DBI} &&  \\
          maximum $q_{c}$ &  0.20   & 0.20    & 0.31    && 0.30   & 0.30    & 0.28  && 0.28  \\ 
          total $q_{c}$   &  45.91  & 47.74   & 87.98   && 86.57  & 82.67   & 75.11 && 69.82\\
        %\hline
                          & \multicolumn{7}{c}{PPI}  &&  \\
         maximum $q_{c}$ & 0.20   & 0.21   & 0.33   && 0.32  & 0.29   & 0.29   && 0.28 \\
         total $q_{c}$   & 47.87  & 50.09  & 97.60  && 92.54 & 81.44  & 78.85  && 69.82 \\ 
        \hline
        \hline
      \end{tabular}                                                                                  
      \caption{Maximum values of $q_{c}$ and the total amount of the cloud mixing ratio at $t = 5$ hours with $nz=600$ points and $\epsilon_{0}=\epsilon_{1}=10^{-5}$.
               The total amount of the cloud mixing ratio is calculated by estimating the integral $q_{c}$.
               The units of $q_{c}$ are $g/kg$.}
      \label{tab:qc1}
    \end{table}

  \section{Discussion and Concluding Remarks}
\label{sec:discussion-conlusion}

  In this work we introduced HiPPIS: a high-order 1D, 2D, and 3D data-bounded and positivity-preserving interpolation software for structured meshes.
  The software implementation is based on the mathematical framework in Section \ref{sec:background} and the algorithms in Section \ref{sec:algorithm}.
  The software is self-contained and can be incorporated into larger codes that require data-bounded or positivity-preserving interpolation.
  The interface is designed to be similar to commonly used PCHIP and splines interfaces.
  The algorithms used in the software extend the DBI and PPI methods introduced in ~\cite{Ouermi2022} by adding more options for the stencil construction process that can be set by the user with the parameter $st$.
  For a given interval $I_{i}$, the algorithm starts with the stencil $\mathcal{V}_{0} = \{x_{i}, x_{i+1}\}$ and successively appends points to the left and/or right of $\mathcal{V}_{0}$ to form the final stencil.
  The stencil construction is done in accordance with the DBI and PPI conditions outlined in Equations (\ref{eq:B_PPIminus}) and (\ref{eq:B_PPIplus}).
  In addition to the different options for the stencil selection process, the software introduces a parameter $\epsilon_{1}$ that can be used to adjust the bounds of the interpolants in the intervals where extrema are detected.
  %
  %\color{red}
  %\st{Strategies for reorganizing the code structure to improve vectorization, locality, and overall computational performance of the DBI and PPI methods are presented in Section} \ref{subsec:implementation-performance}. 
  %\color{black}

  %
  Various 1D and 2D examples are employed to evaluate the use of the DBI and PPI software in different contexts.
  The results in Section \ref{sec:numerical-examples} show that for Examples I, II, IV, and V, which are based on underlying smooth functions, the DBI and PPI methods produce more accurate results compared to PCHIP and MQSI.
  For Example III, and VI which are obtained from discontinuous and $C^{0}$-continuous functions, all the methods have comparable errors.
  Figures \ref{fig:heaviside}-\ref{fig:surface2} show that the parameters $\epsilon_{0}$ and $\epsilon_{1}$ can be used to reduce undesired oscillations from the PPI method.
  We only report the results for $st=3$ because the differences between $st=1$, $st=2$, and $st=3$ are negligible for the examples in Section \ref{sec:numerical-examples}.

  Section \ref{sec:error-mapping-examples} provides an analysis and numerical comparison of the mapping error when the DBI and PPI methods are used to map data values between different meshes.
  The error analysis in Section \ref{subsec:error} shows that it is important to keep mapping errors smaller than the already existing global errors from other calculations.
  Removing negative values and spurious oscillations can help reduce the mapping error.
  
  %
  %Various 1D and 2D examples are employed to evaluate the use of the DBI and PPI software in different contexts.
  %
  The results in Tables \ref{tab:f4-map} and \ref{tab:twp-ice-map} show that using small values for parameters $\epsilon_{0}$ and $\epsilon_{1}$ improves the approximation in cases where the input data are coarse.
  Small values of $\epsilon_{0}$  and $\epsilon_{1}$ further restrict how much the interpolant is allowed to grow beyond the data values.
  The parameters $\epsilon_{0}$ and $\epsilon_{1}$ are used to adjust the lower and upper bounds on each interpolant according to Equations (\ref{eq:umin2}) and (\ref{eq:umax2}). 
  The study of the modified Runge example in Section \ref{subsec:runge} and TWP-ICE example in Section \ref{subsec:twp-ice} demonstrated that for a profile with steep gradients or fronts, more points are required to better take advantage of the DBI and PPI algorithm.
  If there are not sufficiently many points in the regions with steep gradients or fronts, increasing the polynomial degree may not improve the accuracy.
  The results in Tables \ref{tab:f4-map} and \ref{tab:twp-ice-map} show that once the resolution is sufficiently increased, the approximations improve as the polynomial degree increases.

  In the BOMEX test case, prioritizing a symmetry ($st=2$) or locality ($st=3$) leads to better approximations compared to the ENO stencil ($st=1$) using the DBI and PPI methods. 
  Using the ENO stencil ($st=1$) produces significantly less cloud mixing ratio compared to both the prioritizing symmetry and locality.
  In the BOMEX example with parameters $\epsilon_{0}$ and $\epsilon_{1}$ greater than $10^{-5}$, the PPI method allows for oscillations that degrade the approximation compared to the DBI and PCHIP approaches.
  The MQSI method is not used for the BOMEX example because it oversmoothes the different variables at each time step and leads to no production of cloud mixing ratio. 

  In summary, this work provided (1) a high-order DBI and PPI software for 1D, 2D, and 3D structured meshes; (2) an analysis of the mapping error when using the DBI or PPI to map data values between meshes; and (3) an evaluation of the DBI and PPI methods in the context of function approximation and interpolating data values between different meshes. %\color{red}\st{; and (4) code and data restructuring techniques used to enable vectorization, increase locality, and improve overall computational performance.} \color{black}
  As this work continues, we plan to investigate different approaches for extending the DBI and PPI methods to unstructured 2D and 3D meshes.  

\begin{acks}
This work has been supported by the US Naval Research Laboratory (559000669), 
the National Science Foundation (1521748), and the Intel Graphics and Visualization Institute at the University of Utah's Scientific Computing and Imaging (SCI) Institute (29715).
\end{acks}
\vspace{-0.5mm}
%\bibliographystyle{ACM-Reference-Format}
\bibliographystyle{acm}
\bibliography{references}

\end{document}